\documentclass[11pt]{article}
\usepackage{amsfonts}
\usepackage{amssymb}
\usepackage{amsthm}
\usepackage{amsmath}
\usepackage{graphicx}
\usepackage{etoolbox,xstring,mfirstuc,textcase}
\usepackage{empheq}
\usepackage{indentfirst}
\usepackage{cite} 
\usepackage{mathrsfs}
\usepackage{graphics}
\usepackage{xcolor}  
\usepackage{bm}
\usepackage{cases}
\usepackage{times}
\usepackage{changes}
\usepackage{lineno}
\usepackage{hyperref}
\definechangesauthor[name={R1 cusse}, color=orange]{r1}
\definechangesauthor[name={Ty cusse}, color=purple]{ty}

\newcommand{\refform}[1]{\hyperref[#1]{\ref{#1}}}

\newcommand{\norm}[1]{\left\| #1 \right\|}

\newcommand{\hs}[1]{H^{#1}(\Omega)}

\newcommand{\bracket}[1]{\left( #1 \right)}

\newcommand{\limit}[1]{\mathop{\mathrm{lim}}\limits_{#1}}

\newcommand{\abs}[1]{\left| #1 \right|}
\newcommand{\lp}[1]{L^{#1}(\Omega)}

\newtheorem{theorem}{\textbf{Theorem}}[section]
\newtheorem{lemma}{\textbf{Lemma}}[section]
\newtheorem{proposition}{\textbf{Proposition}}[section]
\newtheorem{corollary}{\textbf{Corollary}}[section]
\newtheorem{remark}{\textbf{Remark}}[section]
\newtheorem{definition}{\textbf{Definition}}[section]
\newtheorem{assumption}{\textbf{Assumption}}[section]
\allowdisplaybreaks[4] 
\def\be{\begin{equation}}
\def\ee{\end{equation}}
\def\bea{\begin{eqnarray}}
\def\eea{\end{eqnarray}}
\def\bt{\begin{theorem}}
\def\et{\end{theorem}}
\def\bl{\begin{lemma}}
\def\el{\end{lemma}}
\def\bnum{\begin{numcases}{}}
\def\enum{\end{numcases}}
\def\br{\begin{remark}}
\def\er{\end{remark}}
\def\bp{\begin{proposition}}
\def\ep{\end{proposition}}
\def\bc{\begin{corollary}}
\def\ec{\end{corollary}}
\def\bd{\begin{definition}}
\def\ed{\end{definition}}
\def\ba{\begin{assumption}}
\def\ea{\end{assumption}}
\def\p{\mathrm{\partial}}

\def\half{\frac{1}{2}}

\def\tvp{\widetilde{\varphi}}

\def\uinf{u_{\infty}}
\def\vinf{v_{\infty}}
\begin{document}
	
\title{Regularity and long-time behavior of global weak solutions to a coupled Cahn-Hilliard system: the off-critical case}

\author{
{Bohan Ouyang}
\footnote{School of Mathematical Sciences, Fudan University, Handan Road 220, Shanghai 200433, China.
Email:  \texttt{22110180033@m.fudan.edu.cn}.
}
}

\date{\today}

\maketitle


\begin{abstract}
\noindent We consider a diffuse interface model that describes the macro- and micro-phase separation processes of a polymer mixture. The resulting system consists of a Cahn-Hilliard equation and a Cahn-Hilliard-Oono type equation endowed with the singular Flory-Huggins potential. For the initial boundary value problem in a bounded smooth domain of $\mathbb{R}^d$ ($d\in\{2,3\}$) with homogeneous Neumann boundary conditions for the phase functions as well as chemical potentials, we study the regularity and long-time behavior of global weak solutions in the \emph{off-critical} case, i.e., the mass is not conserved during the micro-phase separation of diblock copolymers. By investigating an auxiliary system with viscous regularizations, we show that every global weak solution regularizes instantaneously for $t>0$. In two dimensions, we obtain the instantaneous strict separation property under a mild growth condition on the first derivative of potential functions near pure phases $\pm 1$, while in three dimensions, we establish the eventual strict separation property for sufficiently large time. Finally, we prove that every global weak solution converges to a single equilibrium as $t\to +\infty$.
 \medskip \\
\noindent
\textbf{Keywords:}  Cahn-Hilliard equation, Cahn-Hilliard-Oono equation, Regularization, Strict separation, Convergence to equilibrium.\medskip \\
\medskip\noindent
\textbf{MSC 2020:} 35A01, 35B40, 35B65, 35K35, 35Q92.
\end{abstract}

\section{Introduction}
\setcounter{equation}{0}
The Cahn-Hilliard equation provides an efficient tool to study the phase separation process of binary mixtures \cite{CH1958,G1996}. It is a representative of the so-called diffuse
interface models that describe the evolution of free interfaces. The diffuse interface approach has the advantages that it avoids the explicit treatment of free interfaces and can handle complex topological changes in a natural way. In the last decades, the Cahn-Hilliard equation and its variants have been successfully applied in many of segregation-driven problems (see \cite{M2019} and the references cited therein).

In the present work, we analyze the following initial-boundary value problem
\begin{align}
     & \partial_t u =\Delta \mu && \mathrm{in} \ \Omega \times (0,T),\label{eq:nvs1}\\
     & \mu = -\epsilon_u^2 \Delta u + \partial_u F(u,v) && \mathrm{in} \ \Omega \times (0,T), \label{eq:nvs2}\\
     & \partial_t v + \sigma (v-c)=\Delta \varphi && \mathrm{in} \ \Omega \times (0,T), \label{eq:nvs3}\\
     & \varphi = -\epsilon_v^2 \Delta v + \partial_v F(u,v) && \mathrm{in} \ \Omega \times (0,T), \label{eq:nvs4}\\
     & \partial_\mathbf{n} u = \partial_\mathbf{n}v = \partial_\mathbf{n} \mu= \partial_\mathbf{n}\varphi= 0 && \mathrm{on} \ \partial \Omega \times (0,T), \label{eq:nvsb} \\
     & (u,v)|_{t=0} = (u_0,v_0) && \mathrm{in} \ \Omega. \label{eq:nvsi}
\end{align}
Here, $\Omega \subset \mathbb{R}^d$ $(d \in\{ 2,3\})$ is a bounded domain with smooth boundary $\partial\Omega$ and $T>0$ is the final time. The vector $\mathbf{n}=\mathbf{n}(x)$ is the unit outer normal vector on $\partial\Omega$ and $\partial_\mathbf{n}$ denotes the outward normal derivative on the boundary.
The coupled system \eqref{eq:nvs1}--\eqref{eq:nvs4} was introduced in \cite{AHTYN2016} to describe the dynamics of a mixture of a homopolymer and a diblock copolymer (with monomers of type A and B) that
undergoes two distinct but simultaneous phase separation processes.
It consists of a Cahn-Hilliard equation and a Cahn-Hilliard-Oono type equation for the phase functions $u$ and $v$, respectively.
The order parameter $u(x,t):\Omega\times[0,T)\to [-1,1]$ denotes the relative fraction difference between the copolymer and the homopolymer, while the order parameter
 $v(x,t):\Omega\times[0,T)\to [-1,1]$ denotes the relative fraction difference between AB components of the diblock copolymer itself. The pure phases $\pm 1$ of $u$ correspond to a homopolymer rich domain $\{u=-1\}$ and a copolymer rich domain $\{u=1\}$. Similarly, the pure phases $\pm 1$ of $v$ correspond to the A-rich domain $\{v=-1\}$ and the B-rich domain $\{v=1\}$. In the diffuse interface framework, both functions $u$, $v$ smoothly transit from $-1$ to $1$ in narrow transition layers, approximating the sharp interfaces with thickness scales $\varepsilon_u, \varepsilon_v>0$, respectively.
 The functions $\mu, \varphi: \Omega\times[0,T)\to \mathbb{R}$ are associated chemical potentials for the macro-separation (between copolymer and homopolymer) and micro-phase separation (between AB blocks) processes. They are given by the variational derivatives of the following free energy functional
 \begin{equation*} 
    \Psi\left(u,v\right) = \int_{\Omega}\left[ \frac{\epsilon_u^2}{2} |\nabla u|^2
    + \frac{\epsilon_v^2}{2} |\nabla v|^2 + F\left(u,v\right) \right]\mathrm{d}x,
\end{equation*}
 where $F$ is a bivariate potential function.
In \cite{AHTYN2016}, the authors considered the following specific form:
$$
F(u,v)= \frac14(u^2-1)^2+ \frac14(v^2-1)^2 + \alpha u v + \beta u^2 v+ \gamma u v^2.
$$
The first two fourth-order polynomials adopt a double-well structure
and have two different minima $\pm 1$ corresponding to the pure phases. The real parameters $\alpha$, $\beta$ and $\gamma$ present the influences of coupling terms in the free energy. They can alter the $(u,v)$-values of the minima of $F(u,v)$, and subsequently affect the confined morphologies of the polymer blend.
It has been shown in \cite{AHTYN2016} that the system \eqref{eq:nvs1}--\eqref{eq:nvs4} is robust enough to predict many kinds of possible morphologies and offers a guideline of
how the system behaves dynamically when the parameters are varied. 
Besides, we mention that the system under consideration is closely related to some diffuse-interface models for binary mixtures with surfactant, see \cite{KK1997,EDA2013,DGW23,Zhu2019} and the references therein.

 Next, let us give some comments on the linear term $\sigma \left(v-c\right)$ in the Cahn-Hilliard-Oono equation \eqref{eq:nvs3}. The parameter $\sigma$ is related
to the bonding between block A and block B in the copolymer such that its value is inversely proportional to the square of the total chain length $N$ (see \cite{CPW09,AHTYN2016}). Here we treat a general case with $c\in (-1,1)$ being a prescribed constant (see \cite{DG2022,GGM2017}).
 For any given (sufficiently regular) function $g$, denote its spatial mean value by $\overline{g}= |\Omega|^{-1}\int_\Omega g\,\mathrm{d}x$.
  With this notation, we find that the solution $(u,v)$ to problem \eqref{eq:nvs1}--\eqref{eq:nvsi} formally satisfies the ordinary differential equations
\begin{equation*}
\frac{\mathrm{d} \overline{u}}{\mathrm{d}t}=0\quad \text{and}
\quad \frac{\mathrm{d} \overline{v}}{\mathrm{d}t}
+ \sigma \left(\overline{v}-c\right) =0,\quad \forall\, t\in (0,T),
\end{equation*}
with initial values $\overline{u}(0)=\overline{u_0}$, $\overline{v}(0)=\overline{v_0}$.
Then we have
\begin{equation}\label{eq:valuesofaveuv}
\overline{u}(t) = \overline{u_0}\quad \text{and}\quad
\overline{v}(t) = \overline{v_0} e^{-\sigma t}+ c\left(1-e^{-\sigma t}\right),
\quad \forall\, t\in [0,T).
\end{equation}
Hence, $\overline{v}$ will remain constant provided that  $\sigma=0$ or $\overline{v_0} = c$ (when $\sigma\neq 0$). This is usually referred to as the \emph{conserved} case \cite{GGM2017}. On the other hand, the case with $\overline{v_0} \not = c$ is called the \emph{off-critical} case. In particular, if $\sigma>0$, then $\overline{v}(t)$ converges exponentially fast to $c$ as $t\to +\infty$ according to \eqref{eq:valuesofaveuv}.
It is easy to verify that if $c,\overline{u}_0,\overline{v}_0 \in \left(-1+m,1-m\right)$ for some given constant $m \in \left(0,1\right)$, then in both cases the mean values $\overline{u},\overline{v}$ will stay in $\left(-1+m,1-m\right)$ on the whole interval of existence for $(u,v)$.
In \cite{AHTYN2016}, the authors considered the conserved case with $c=\overline{v_0}$. Thanks to \eqref{eq:valuesofaveuv}, the equation \eqref{eq:nvs3} reduces to
$$
 \frac{\mathrm{\partial} v}{\mathrm{\partial} t} + \sigma \left(v-\overline{v}\right)=\Delta \varphi.
$$
Thus, the coupled system \eqref{eq:nvs1}--\eqref{eq:nvs4} can be viewed as a
gradient flow of the following free energy with a nonlocal term:
 \begin{equation}
 \widetilde{\Psi}(u,v)
 = \int_{\Omega}\left[ \frac{\epsilon_u^2}{2} |\nabla u|^2
    + \frac{\epsilon_v^2}{2} |\nabla v|^2 + F(u,v)+ \frac{\sigma}{2} |(-\Delta)^{-\frac12}(v-\overline{v})|^2 \right]\mathrm{d}x,
    \label{OK}
\end{equation}
 where $\Delta$ is the Laplace operator with homogeneous Neumann boundary condition in our current setting (cf. \eqref{eq:nvsb}). Moreover, the system satisfies the energy dissipative law in the conserved case (at least formally):
 \begin{equation}
 \frac{\mathrm{d}}{\mathrm{d} t}\widetilde{\Psi}(u,v) + \int_\Omega \left(|\nabla \mu|^2 + |\nabla \widetilde{\varphi}|^2 \right) \mathrm{d}x =0,
 \label{con-diss}
 \end{equation}
 with $\widetilde{\varphi}=\varphi+ \sigma(-\Delta)^{-1}(v-\overline{v})$.
 The modified free energy \eqref{OK} is closely related to the Ohta-Kawasaki functional for diblock copolymers \cite{CPW09,NO95,OK86}. When $\sigma\neq 0$, the Oono's term $\sigma (v-\overline{v})$ yields possible long range interactions that can generate a variety of minimizers with fine structure in the micro-phase separation process, such as layers, onions and multipods (see \cite{AHTYN2016} for details). For mathematical analysis of the Cahn-Hilliard-Oono equation and its variants, we refer to
 \cite{CGRS2022,FS2024,GGM2017,H2022,M2011} and the references therein.

In this study, we are interested in the theoretical analysis of problem \eqref{eq:nvs1}--\eqref{eq:nvsi}. It is well-known that as a fourth-order parabolic equation, the Cahn-Hilliard (or Cahn-Hilliard-Oono) equation with a regular potential (e.g., a polynomial) does not maintain the maximum principle, that is, its solution may not stay in the physically relevant interval $[-1, 1]$ throughout the evolution, see \cite{M2019} and the references therein. On the other hand, the polynomial double-well potentials used in \cite{AHTYN2016} are just convenient approximations of the Flory-Huggins free energy density for polymers \cite{Flo42,Hug41}:
\begin{align}
S_{\textrm{FH}}(r)=\frac{\theta}{2} \left[(1+r)\ln(1+r) + (1-r)\ln(1-r)\right]- \frac{\theta_{0}}{2} r^2,\quad \forall\, r\in (-1,1).
\label{FH}
\end{align}
The two parameters $\theta$, $\theta_0$ in \eqref{FH} denote the absolute temperature and the critical temperature for phase separation, respectively.
If $0<\theta<\theta_0$, the potential $S_{\textrm{FH}}$ has a double-well structure with two minima inside $(-1,1)$. The singular nature of $S_{\textrm{FH}}$ near the endpoints $\pm 1$ (i.e., pure phases) can ensure the existence of physical solutions with values in $[-1,1]$ (cf. \cite{AW2007,GGM2017}). Inspired by this, the authors of \cite{DG2022} analyzed problem \eqref{eq:nvs1}--\eqref{eq:nvsi} with the choice
\begin{align}
    &F(u,v) = S\left(u;\theta_u,\theta_{0,u}\right) + S\left(v;\theta_v,\theta_{0,v}\right) + W(u,v),
    \notag
\end{align}
where $W(u,v)=\alpha u v + \beta u^2 v+ \gamma u v^2$ and
\begin{align}
    &S\left(r;\theta_r,\theta_{0,r}\right) = \frac{\theta_r}{2} \left[\left(1+r\right)\mathrm{ln}\left(1+r\right) + \left(1-r\right)\mathrm{ln}\left(1-r\right)\right] - \frac{\theta_{0,r}}{2} r^2,
    \notag
\end{align}
with $0 < \theta_r < \theta_{0,r}$, $r\in\{u, v\}$.
For both conserved and off-critical cases, they proved that the resulting initial boundary value problem admits a unique global weak solution $(u,v)$ on $[0,+\infty)$ with certain dissipative estimates in two and three dimensions.
Then they focused on analysis of the conserved case and showed that every weak solution regularizes instantaneously for $t>0$.
When the spatial dimension is two, they further established the instantaneous strict separation property for the macro- and micro-phases, that is,
for any $t_0 > 0$, there exist constants $\omega_u, \omega_v\in (0, 1)$ depending on $t_0$ and the initial data such that
$$
\|u(t)\|_{L^\infty(\Omega)}\leq 1-\omega_u,\quad \|v(t)\|_{L^\infty(\Omega)}\leq 1-\omega_v,\quad \forall\, t\geq t_0.
$$
The above regularization properties allowed them to study the long-time behavior of global weak solutions in both two and three dimensions, i.e., the convergence to a single equilibrium as $t\to +\infty$.

As pointed out in \cite[Section 5]{DG2022}, the nonlinear coupling term $W(u,v)$ in the mixing entropy leads to extra difficulties in the study of regularity properties of global weak solutions (and thus their long-time behavior). Because of the lack of control
on $L^\infty$-norms for approximating phase functions, higher-order estimates of their  time derivatives cannot be achieved by exploiting the Galerkin approximation scheme like in \cite{GGM2017}, where the single Cahn-Hilliard-Oono equation was analyzed.
Inspired by \cite{GGW2018}, the authors of \cite{DG2022} derived higher-order time regularity of weak solutions in the conserved case, using the difference quotients in time.
Unfortunately, this argument seems not valid in the off-critical case.
Comparing with the conserved case, one of the main difficulties came from an additional inner product term $\big(\partial_t^h \varphi , \overline{\partial_t^h v} \big)$, which was due to the possible mass change in the micro-phase separation (see \eqref{eq:valuesofaveuv}). Here, for any function $f:[0,T]\to X$, with $X$ being a real Banach space, we denote $\partial_t^h f= h^{-1}[f(t+h)-f(t)]$ for any $h > 0$ and $t\geq 0$. Due to the difficulty mentioned above, regularity and long-time behavior  of global weak solutions to problem \eqref{eq:nvs1}--\eqref{eq:nvsi} were only analyzed in the conserved case, that is, the masses are conserved and the energy is dissipative (see \eqref{con-diss}). 

Our aim in the present contribution is to study the regularity properties and long-time behavior of global weak solutions to problem \eqref{eq:nvs1}--\eqref{eq:nvsi} in the off-critical case.

Under some general assumption on the nonlinearities, we recover the existence of global weak solutions (see Proposition \ref{cdivnv}) and show that every global weak solution regularizes instantaneously for $t>0$ (see Theorem \ref{hoenv}).
Our proof is inspired by \cite{MZ2004}, that is, we investigate an auxiliary system with viscous regularization in the chemical potentials $\mu$ and $\varphi$. Given sufficiently regular initial data being strictly separated from $\pm 1$, the coupled viscous Cahn-Hilliard system possesses global strong solutions that are smooth enough for us to perform higher-order estimates (see Proposition \ref{ssve}). In particular, thanks to the viscous regularizing terms, we can apply the method in \cite{MZ2004} to obtain the instantaneous strict separation property for phase functions $u, v$ in both two and three dimensions. Combining uniform estimates for the approximate solutions with a compactness argument, we are able to prove the existence of global weak solutions to the original problem \eqref{eq:nvs1}--\eqref{eq:nvsi} and their instantaneous regularization property. In the two dimensional case, we further establish the instantaneous strict separation of weak solutions (see Theorem \ref{strictseparation2d}),  by applying a De Giorgi type iteration scheme proposed in the recent work \cite{GP2023}.
With the aid of this approach, we achieve the result under weaker assumptions on the singular potentials than those in \cite{DG2022}. Next, taking advantage of the viscous regularization for the system and the above mentioned regularizing effects, we are able to characterize the $\omega$-limit set (see Proposition \ref{omegalim}) and derive the eventual strict property of global weak solutions in both two and three dimensions (see Theorem \ref{longtimesep}). It is worth mentioning that the validity of this property in the off-critical case was open even for the single Cahn-Hilliard-Oono equation (see \cite{GGM2017}). Finally, we derive an extended {\L}ojasiewicz-Simon inequality that works in the situation with possible mass change (see Proposition \ref{ls}). This extends the corresponding result in \cite[Proposition 7.2]{DG2022} for the conserved case and enables us to establish the convergence to a single equilibrium (see Theorem \ref{convergequi}).

\textit{Plan of this paper.} In Section \ref{mr}, we first introduce the functional setting some notations, then we state the main results of this study. In Section \ref{ecvs} we analyze the viscous regularizing system and derive uniform estimates that are independent of the approximating parameter. In Section \ref{nv}, we prove the existence and regularity of global weak solutions to the original problem. Besides, we establish the instantaneous  strict separation property in the two dimensional case. 
In Section \ref{longt}, we study the $\omega$-limit set and show that every weak solution becomes strictly separated for sufficiently large time. Then we prove the convergence to an equilibrium as $t\to +\infty$ with the aid of the {\L}ojasiewicz-Simon inequality. In the appendix, we sketch a proof for the well-posedness of an auxiliary problem with regular potentials and present some useful results on an nonlinear elliptic Neumann problem with singular potential. 


\section{Main Results}\label{mr}
\setcounter{equation}{0}
\subsection{Preliminaries} 
We first introduce some notations and conventions.
Let $\mathcal{X}$ be a (real) Banach space with norm $\left\|\cdot\right\|_{\mathcal{X}}$. We use $\mathcal{X}^*$, $\langle \cdot,\cdot \rangle_{\mathcal{X}^*,\mathcal{X}}$ to represent its dual space
and the associated duality pairing. For a Hilbert space $\mathcal{H}$, we denote the associated inner product by $(\cdot,\cdot)_{\mathcal{H}}$.
Throughout this paper, we assume that $\Omega \subset \mathbb{R}^d$ ($d\in\{2,3\}$) is a bounded domain with smooth boundary $\partial \Omega$.
For the standard Lebesgue and Sobolev spaces on $\Omega$, we use the notations $L^{p}(\Omega)$, $W^{k,p}(\Omega)$ for any $p \in [1,+\infty]$ and $k\in \mathbb{N}$,
equipped with the corresponding norms
$\|\cdot\|_{L^{p}(\Omega)}$, $\|\cdot\|_{W^{k,p}(\Omega)}$, respectively. When $p = 2$, these spaces are Hilbert spaces and we use the standard convention $H^{k}(\Omega) := W^{k,2}(\Omega)$. For simplicity, the norm and inner product of $L^2(\Omega)$ will be denoted by $\|\cdot\|$ and $(\cdot,\cdot)$, while the pairing between $H^1(\Omega)$ and $H^1(\Omega)^*$ will be denoted by $\langle \cdot,\cdot \rangle$.
Bold letters will be used for vector-valued
spaces, for instance, we denote the vector-valued Lebesgue spaces by $\bm{L}^p(\Omega)= L^p(\Omega;\mathbb{R}^d)$, $p \in [1,+\infty]$.

Given a measurable set $I$ of $\mathbb{R}$, we introduce the function space $L^p(I;\mathcal{X})$ with $p\in [1,+\infty]$, which consists of Bochner measurable $p$-integrable functions (if $p \in [1,+\infty)$) or essentially bounded functions (if $p =+\infty$) with values in a given Banach space $\mathcal{X}$. If $I=(a,b)$, we write for simplicity $L^p(a,b;\mathcal{X})$. The space $L^p_{\mathrm{uloc}}(0,+\infty;\mathcal{X})$ denotes the uniformly local variant of $L^p(0,+\infty;\mathcal{X})$ consisting of all strongly measurable $f:[0,+\infty)\to \mathcal{X}$ such that
\begin{equation*}
    \|f\|_{L^p_{\mathrm{uloc}}(0,+\infty);\mathcal{X})}:= \mathop\mathrm{sup}\limits_{t\ge0} \|f\|_{L^p(t,t+1;\mathcal{X})} < \infty.
\end{equation*}
If $T\in(0,+\infty)$, we find $L^p_{\mathrm{uloc}}(0,T;\mathcal{X})=L^p(0,T;\mathcal{X})$.

The following shorthands will be frequently used
$$
H := L^2\left(\Omega\right),
\quad V := H^1\left(\Omega\right),
\quad W :=\big\{u\in H^2(\Omega)\ |\ \partial_\mathbf{n}u=0\ \text{on}\ \partial\Omega\big\}.
$$
As usual, $H$ is identified with its dual. We have the continuous, dense, and compact embeddings:
\begin{equation*}
    W \hookrightarrow V  \hookrightarrow H \hookrightarrow V^*.
\end{equation*}
Besides, the interpolation inequality holds (see \cite{GGM2017})
\begin{equation*} 
    \left\|f\right\|^2 \le \xi \left\|\nabla f\right\|^2 + C\left(\xi\right) \left\|f\right\|^2_{V^*}, \quad \forall f \in V,
\end{equation*}
where  $\xi \in (0,1)$ is arbitrary and $C\left(\xi\right)$ is a positive constant only depending on $\xi$ and $\Omega$.
For every $f\in V^*$, we denote by $\overline{f}$ its generalized mean value over $\Omega$ such that
$\overline{f}=|\Omega|^{-1}\langle f,1\rangle$; if $f\in L^1(\Omega)$, its mean value is simply given by $\overline{f}=|\Omega|^{-1}\int_\Omega f \,\mathrm{d}x$.
Then we recall the well-known Poincar\'{e}-Wirtinger inequality:
\begin{equation*}
\left\|f-\overline{f}\right\|\leq C_P \|\nabla f\|,\quad \forall\,
f\in V,
\end{equation*}
where the positive constant $C_P$ depends only on $\Omega$.

We denote by $A_N: V \to V^*$ the extension of the minus Laplace operator subject to the homogeneous Neumann boundary condition such that
\begin{equation*}
    \left\langle A_N f,g \right\rangle = \int_{\Omega} \nabla f \cdot \nabla g \,\mathrm{d}x, \quad \forall\,f,g \in V.
\end{equation*}
Set the linear subspaces
\begin{equation*}
 H_{0}:=\{f\in H\ |\ \overline{f} =0\},\quad
 V_0 := V\cap H_0,\quad
 V^*_0 := \left\{ L \in V^* \ |\   \left\langle L,1 \right\rangle = 0 \right\}.
\end{equation*}
It is easy to verify that the restriction of $A_N$ to $V_0$ is an isomorphism between $V_0$ and $V^*_0$. Thus, we can define the inverse operator $\mathcal{N} := \left(A_N|_{V_0}\right)^{-1}:V_0^* \to V_0$. The following identities hold (see \cite{DG2022,MZ2004})
\begin{align*}
&\left\langle A_Nu,\mathcal{N}L \right\rangle = \left\langle L,u \right\rangle,\quad \forall\, u \in V_0,\ L \in V_0^*,
\\
&\left\langle L_1,\mathcal{N}L_2 \right\rangle = \left\langle \nabla \left(\mathcal{N}L_1\right), \nabla \left(\mathcal{N}L_2\right) \right\rangle,
\quad \forall\, L_1,\,L_2\ \in V_0^*.
\end{align*}
Define
\begin{align*}
    &\left\|L\right\|_* := \left\| \nabla \left(\mathcal{N}L\right)\right\| = \sqrt{\left\langle L,\mathcal{N}L \right\rangle},\quad \forall\, L\in V^*_0, \\
    &\left\|L\right\|_{-1}^2 := \left\|L-\overline{L}\right\|_*^2 + |\overline{L}|^2, \quad \forall\, L\in V^*.
\end{align*}
We find that $\left\|\cdot\right\|_*$ and $\left\|\cdot\right\|_{-1}$ are equivalent norms on $V^*_0$ and $V^*$ with respect to the usual dual norms (see \cite{MZ2004}).

In the subsequent analysis, the capital letter $C$ will denote a generic positive constant that depends on the structural data of the problem. Its meaning may change from line to line and even within the same chain of computations. Specific dependence will be pointed out if necessary.

\subsection{Statement of results}

Let us present the assumptions for problem \eqref{eq:nvs1}--\eqref{eq:nvsi}.
We recall that the parameters $\varepsilon_u, \varepsilon_v>0$ are proportional to the thickness of transition layers between different components in the polymer mixture, and thus are related to the rapidity of variation of $u$ and $v$ in the interfacial region (see \cite{AHTYN2016}). Nevertheless, their values will not affect the subsequent analysis in this study, since we always work with fixed positive $\varepsilon_u, \varepsilon_v$ and do not consider the asymptotic behavior as $\varepsilon_u, \varepsilon_v\to 0$ (i.e., the sharp-interface limit).
Hence, without loss of generality, we make the following assumption on parameters of the system:
\begin{itemize}
\item[$(\mathbf{H0})$]
$\varepsilon_u = \varepsilon_v = 1$,\ $\sigma>0$,\ $c\in (-1,1)$, \ $0<\theta_u<\theta_{0,u}$, \ $0<\theta_v<\theta_{0,v}$.
\end{itemize}
In the subsequent analysis, we consider the free energy 
\begin{equation}\label{eq:energyfunctionalr}
    \Psi(u,v) = \int_{\Omega}\left[ \frac{1}{2} |\nabla u|^2
    + \frac{1}{2} |\nabla v|^2 + F(u,v) \right] \mathrm{d}x,
\end{equation}
with the bivariate potential
\begin{align}
    &F(u,v) = S\left(u;\theta_u,\theta_{0,u}\right) + S\left(v;\theta_v,\theta_{0,v}\right) + W(u,v),
    \label{eq:F}
\end{align}
where
\begin{align*}
    &S\left(u;\theta_u,\theta_{0,u}\right) =  \widehat{S}_{(u)}(u) -\frac{\theta_{0,u}}{2} u^2,
    \quad S\left(v;\theta_v,\theta_{0,v}\right) =  \widehat{S}_{(v)}(v)- \frac{\theta_{0,v}}{2} v^2.
\end{align*}
In view of \cite{DG2022,GP2023}, we impose the following structural assumptions on the nonlinear functions $\widehat{S}_{(j)}$ $(j\in\{u,v\})$ and $W$:
\begin{enumerate}
    \item[$(\mathbf{H1})$] For $j\in\{u,v\}$, $\widehat{S}_{(j)}\in C\left([-1,1]\right)\cap C^2(-1,1)$, and \begin{align*}
        & \lim_{s \to -1^+} \widehat{S}'_{(j)}(s) = -\infty,\qquad \lim_{s \to 1^-} \widehat{S}'_{(j)}(s) = +\infty,\\
        & \widehat{S}_{(j)}''(s) \ge \theta_j>0, \quad \forall\,s\in (-1,1).
        \end{align*}
    We extend $\widehat{S}_{(j)}(s)=+\infty$ for all $|s|>1$ and without loss of generality, we assume  $\widehat{S}_{(j)}(0)=\widehat{S}'_{(j)}(0)=0$.
    \item[$(\mathbf{H2})$] There exists some constant $\rho> 1/2$ such that as $\delta \to 0^+$, it holds
    \begin{equation}\notag 
        \frac{1}{\widehat{S}_{(j)}'(1-2\delta)} = O\left(\frac{1}{| \ln \delta|^{\rho}}\right),
        \quad \frac{1}{\widehat{S}_{(j)}'(-1+2\delta)} = O\left(\frac{1}{| \ln \delta|^{\rho}}\right).
    \end{equation}
    \item[$(\mathbf{H3})$] $W\in C^2(\mathbb{R}^2;\mathbb{R})$.
\end{enumerate}

\begin{remark}\rm
    The logarithmic potential (cf. \eqref{FH})
    \begin{equation}\notag 
        \widehat{S}_{(j)}(s) = \frac{\theta_j}{2} \left[\left(1+s\right)\ln(1+s) + (1-s)\ln(1-s)\right], \quad \forall\,s \in (-1,1), \ j\in \{u,v\},
    \end{equation}
    fulfills the assumptions $(\mathbf{H1})$--$(\mathbf{H2})$ with $\rho=1$, while the bivariate polynomial $W(u,v)=\alpha u v + \beta u^2 v+ \gamma u v^2$ considered in \cite{AHTYN2016,DG2022} fulfills the assumption  $(\mathbf{H3})$. The assumption $\widehat{S}_{(j)}(0)=\widehat{S}'_{(j)}(0)=0$ and the convexity also imply $\widehat{S}_{(j)}\geq 0$ on $[-1,1]$. The assumption $(\mathbf{H2})$ characterizes the growth of the first order derivative of singular potentials $\widehat{S}_{(j)}$ and it only plays a role in the proof for the instantaneous strict separation property of global weak solutions in two dimensions (see Proposition \ref{strictseparation2d}).
    It is worth mentioning that our assumptions on the nonlinearities are weaker than those in \cite{DG2022}, cf. \cite[Section 4.1]{DG2022} and also \cite[Assumption A]{DG2022}.
\end{remark}

Next, we introduce the definition of weak solutions to problem \eqref{eq:nvs1}--\eqref{eq:nvsi} (cf.  \cite[Definition 3.1]{DG2022}).
\bd\label{wsnvd}
Let $\Omega \subset \mathbb{R}^d$ $(d\in\{2,3\})$ be a bounded domain with smooth boundary $\partial \Omega$ and $T \in \left(0,+\infty\right]$. Assume that $(\mathbf{H0})$, $(\mathbf{H1})$ and $(\mathbf{H3})$ are satisfied. For any initial data $u_0, v_0 \in V$ such that $F\left(u_0,v_0\right) \in L^1\left(\Omega\right)$ and $\overline{u_0},\overline{v_0} \in \left(-1,1\right)$, the quadruple $(u,v,\mu,\varphi)$ is called a weak solution to problem \eqref{eq:nvs1}--\eqref{eq:nvsi} on $[0,T]$, if the  following properties are satisfied:
\begin{align*}
& u, v \in L^{\infty}(0,T;V)\cap C([0,T];H)
\cap L^2_{\mathrm{uloc}}(0,T;W) \cap H^1_{\mathrm{uloc}}(0,T;V^*),\\
& u, v \in L^\infty(\Omega\times (0,T))\ \ \text{with}\ |u(x,t)|, |v(x,t)|<1\ \ \text{a.e. in}\ \Omega\times (0,T),\\
&\mu, \varphi \in  L^2_{\mathrm{uloc}}(0,T;V),\quad \widehat{S}'_{(u)}(u),\ \widehat{S}'_{(v)}(v)  \in L^2_{\mathrm{uloc}}(0,\infty;H),
\end{align*}
and
\begin{align*}
    & \left\langle \partial_t u, \eta \right\rangle + (\nabla \mu, \nabla\eta) = 0,&&\forall\,\eta\in V, \ \ \text{a.e. in}\ (0,T),\\
    & \left\langle \partial_t v, \eta \right\rangle + \sigma(v-c,\eta) + (\nabla \varphi, \nabla \eta) = 0,&&\forall\,\eta\in V, \ \ \text{a.e. in}\ (0,T),\\
    &\mu= -\Delta u+ \partial_u F(u,v), &&
     \text{a.e. in}\ \Omega\times (0,T), \\
    &\varphi= -\Delta v+ \partial_v F(u,v),&&
     \text{a.e. in}\ \Omega\times (0,T),
\end{align*}
moreover,
$\partial_\mathbf{n} u = \partial_\mathbf{n} v= 0$ almost every on $\partial \Omega\times (0,T)$,
and $(u(0),v(0)) = (u_0,v_0)$ almost every in $\Omega$.
\ed

Then we have the following result on the existence and uniqueness of a global weak solution:
\bp \label{cdivnv}
Let $\Omega \subset \mathbb{R}^d$ $(d\in\{2,3\})$ be a bounded domain with smooth boundary $\partial \Omega$. Assume that $(\mathbf{H0})$, $(\mathbf{H1})$ and $(\mathbf{H3})$ are satisfied.

\begin{itemize}
\item[(1)] For any initial data $u_0, v_0 \in V$ such that $F(u_0,v_0) \in L^1(\Omega)$ and $\overline{u_0},\overline{v_0} \in (-1,1)$,
problem \eqref{eq:nvs1}--\eqref{eq:nvsi} admits a unique global weak solution  $(u,v,\mu,\varphi)$ on $[0,+\infty)$ in the sense of Definition \ref{wsnvd}.
Moreover, the following energy inequality holds
\begin{align}\label{eq:denvs-1}
     &   \widehat{\Psi}(u(t),v(t))  \le \widehat{\Psi}(u_0,v_0) e^{-(1+\sigma) t} + C_1,\quad \forall\, t\ge0, \\
     &   \frac{1}{2} \int_t^{t+1} \left( \| \nabla \mu (\tau ) \|^2 + \| \nabla \varphi (\tau) \|^2  \right) \mathrm{d}\tau \leq \widehat{\Psi}(u_0,v_0) e^{-(1+\sigma) t} + C_2,\quad \forall\, t\ge0, \label{eq:denvs-2}
\end{align}
where
\begin{align}
        \widehat{\Psi}(u,v) & = \Psi(u,v) + \frac{1}{2} \left(\| u-\overline{u}\|_{V_0^*}^2 +
        \| v-\overline{v} \|_{V_0^*}^2 \right),
        \notag
    \end{align}
with $\Psi(u,v)$ given by \eqref{eq:energyfunctionalr}, $C_1$, $C_2$ are positive constants that only depend on parameters of the system, $\Omega$ and $\overline{u_0}$, $\overline{v_0}$.
\item[(2)] Given $R\ge 0$, $T\in (0,+\infty)$ and $m \in (|c|,1)$, there exists a positive constant $C$ depending on $m$, $R$, $T$ such that, for any solutions $(u_1,v_1)$, $(u_2,v_2)$ on $[0,T]$ originating from the initial data $(u_{01},v_{01})$, $(u_{02},v_{02})$ satisfying $\Psi(u_{0i},v_{0i})\le R$ and $|\overline{u_{0i}}|$, $|\overline{v_{0i}}|\le m$ $(i=1,2)$, the continuous dependence estimate
\begin{equation}
    \begin{aligned}
        & \| u_1(t) - u_2(t) \|_{V^*}^2 + \| v_1(t) - v_2(t) \|_{V^*}^2
        +  \int_0^T \left(\| u_1(t) - u_2(t) \|_{V}^2+ \| v_1(t) - v_2(t)\|_{V}^2\right) \mathrm{d}t \\
        &\quad \le  C \left(\| u_{01} - u_{01} \|_{V^*}^2 + \| v_{01} - v_{02}\|_{V^*}^2 + |\overline{u_{01}}-\overline{u_{02}}| + |\overline{v_{01}}-\overline{v_{02}}|\right)
        \label{es-conti1}
    \end{aligned}
\end{equation}
holds for every $t\in \left[0,T\right]$.
\end{itemize}
\ep

\begin{remark}\rm
In \cite[Theorem 3.1]{DG2022}, the authors have proved a global well-posedness result similar to Proposition \ref{cdivnv} for $\sigma\geq 0$, under some assumptions that were slightly stronger than $(\mathbf{H1})$ and $(\mathbf{H3})$. In our current setting, we recover the existence result by a different approach. Since the continuous dependence estimate \eqref{es-conti1} can be obtained by exactly the same argument as in \cite[Section 4.4]{DG2022}, its proof will be omitted.
Although our focus in this paper is the off-critical case with $\sigma>0$, the results in Proposition \ref{cdivnv} can be easily extended to the conserved case as well.
\end{remark}

\begin{remark}\rm
    Following \cite[Remark 3.3]{GGG2017}, we find that $t \to \|u(t)\|_{L^{\infty}(\Omega)}$ is measurable and essentially bounded, and the same conclusion holds for $v$ (cf. \cite[Remark 3.2]{DG2022}).
    Arguing as in \cite{A2009,GGW2018,GMT2019}, we also obtain $u,v \in L^4_{\mathrm{uloc}}(0,+\infty;W)$. Furthermore, an application of Lemma \ref{nonlinearelliptic} yields
    \begin{equation*}
        u,v \in L^2_{\mathrm{uloc}}(0,+\infty;W^{2,p}(\Omega)), \quad \widehat{S}'_{(u)}(u),\ \widehat{S}'_{(v)}(v)  \in L^2_{\mathrm{uloc}}(0,+\infty;L^p(\Omega)),
    \end{equation*}
    where $p = 6$ if $d=3$, or $p \in [2,+\infty)$ if $d=2$.
\end{remark}

Now we are in a position to state our main results.

\bt[Instantaneous regularization]
\label{hoenv}
Suppose that the assumptions of Proposition \ref{cdivnv} are satisfied.
Let $(u, v, \mu, \varphi)$ be the unique global weak solution obtained therein.
Then for any  $\kappa > 0$, we have
\begin{align*}
    & u,v \in L^{\infty}( \kappa,+\infty; W^{2,p}(\Omega)),\quad
      \partial_t u, \partial_t v \in L^{\infty}( \kappa,+\infty;V^*) \cap L^2_{\mathrm{uloc}}( \kappa,+\infty;V ),\\
    & \mu,\varphi \in L^{\infty}( \kappa,+\infty;V )\cap L^2_{\mathrm{uloc}}( \kappa,+\infty;H^3(\Omega)),
\end{align*}
where $p = 6$ if $d=3$, or $p \in [2,+\infty)$ if $d=2$.
Thus, the equations \eqref{eq:nvs1}, \eqref{eq:nvs3} are satisfied
almost everywhere in $\Omega\times [\kappa,+\infty)$ and the boundary conditions $\partial_\mathbf{n} \mu=\partial_\mathbf{n} \varphi= 0$ hold
almost everywhere on $\p \Omega\times [\kappa,+\infty)$.
Moreover, it holds
\begin{align}
   & \|u\|_{L^\infty(\kappa,t; W^{2,p}(\Omega))}
   + \|v\|_{L^\infty(\kappa,t; W^{2,p}(\Omega))}
   + \|\mu \|_{L^\infty(\kappa,t;V)}
   + \|\varphi \|_{L^\infty(\kappa,t;V)} \notag \\
   &\quad
   + \int_t^{t+1} \left(\| \partial_t u(\tau) \|_V^2 + \| \partial_t v(\tau) \|_V^2
   +\|\mu(\tau)\|_{H^3(\Omega)}^2+\|\varphi(\tau)\|_{H^3(\Omega)}^2 \right) \mathrm{d}\tau
   \le C,\quad \forall\, t \ge \kappa,
   \notag
   \end{align}
   where the positive constant $C$ depends on ${\Psi}(u_0,v_0)$, $\overline{u_0}$, $\overline{v_0}$, $\Omega$, parameters of the system and $\kappa$.
\et

\bt[Instantaneous strict separation in two dimensions]\label{strictseparation2d}
Suppose that the assumptions of Proposition \ref{cdivnv} are satisfied. In addition, we assume that $d=2$ and $(\mathbf{H2})$ is satisfied. Let $(u, v)$ be the unique global weak solution to problem \eqref{eq:nvs1}--\eqref{eq:nvsi}. Then for any $\kappa \in(0,1]$, there exists $\delta_\kappa \in (0,1)$ such that
\begin{equation}
    \| u(t) \|_{C(\overline{\Omega})} \le 1 - \delta_\kappa,
    \quad \|v(t)\|_{C(\overline{\Omega})} \le 1 - \delta_\kappa,
    \quad \forall\,t \ge \kappa,
    \label{isp-2d}
\end{equation}
where the constant $\delta_\kappa$ depends on ${\Psi}(u_0,v_0)$, $\overline{u_0}$, $\overline{v_0}$, $\Omega$, parameters of the system and $\kappa$.
\et

\bt[Eventual strict separation in two and three dimensions]\label{longtimesep}
Suppose that the assumptions of Proposition \ref{cdivnv} are satisfied.
Let $(u, v)$ be the unique global weak solution to problem \eqref{eq:nvs1}--\eqref{eq:nvsi}. There exist $\delta_{\mathrm{SP}} \in (0,1)$ and $T_{\mathrm{SP}}\gg 1$ such that
\begin{equation}
     \| u(t) \|_{C(\overline{\Omega})} \le 1 - \delta_{\mathrm{SP}},
    \quad \|v(t)\|_{C(\overline{\Omega})} \le 1 - \delta_{\mathrm{SP}},
    \quad \forall\,t \ge T_{\mathrm{SP}}.
    \label{esp-3d}
\end{equation}
\et

Thanks to the above regularizing properties, we are able to address the long-time behavior of global weak solutions under the following additional assumption:
\begin{itemize}
\item[$\mathbf{(H4)}$] $\widehat{S}_{(u)}$, $\widehat{S}_{(v)}$ are real analytic on $(-1,1)$ and $W$ is real analytic on $(-1,1)^2$.
\end{itemize}

\bt[Convergence to equilibrium]\label{convergequi}
Suppose that the assumptions of Proposition \ref{cdivnv} are satisfied. Let $(u, v)$ be the unique global weak solution to problem \eqref{eq:nvs1}--\eqref{eq:nvsi}. If in addition, $\mathbf{(H4)}$ is fulfilled.
Then we have
\begin{equation*}
    \limit{t \to +\infty} \norm{(u(t),v(t))-(\uinf,\vinf)}_{\hs{2-\epsilon} \times \hs{2-\epsilon}}  = 0,
\end{equation*}
for any $\epsilon \in (0,1/2)$. Here, $(\uinf,\vinf)$ is a steady state that satisfies
 \begin{equation}
 \begin{cases}
     -\Delta \uinf + \partial_u F(\uinf,\vinf) = \overline{\partial_u F(\uinf,\vinf)},\qquad \qquad \qquad\qquad\quad  \ \, \mathrm{in} \ \Omega, \\
     -\Delta \vinf + \partial_v F(\uinf,\vinf) +\sigma \mathcal{N}\left(\vinf-\overline{\vinf}\right) = \overline{\partial_v F(\uinf,\vinf)},\quad\quad    \mathrm{in} \ \Omega, \\
     \partial_\mathbf{n}\uinf = \partial_\mathbf{n}\vinf = 0,
      \qquad \qquad \qquad\qquad \qquad \qquad\qquad \qquad \quad\quad  \ \ \mathrm{on} \ \p \Omega, \\
      \text{with mass constraints}\quad \overline{\uinf}=\overline{u_0},\quad \overline{\vinf}=c.
\end{cases}
\label{sta0}
\end{equation}
\et
\begin{remark}\rm 
By the same argument as in \cite{DG2022} (see also \cite{JWZ2015}), we can derive an estimate on the convergence rate:
$$
\norm{(u(t),v(t))-(\uinf,\vinf)}_{V^*\times V^*}\leq C(1+t)^{-\frac{\theta}{1-2\theta}},\quad \forall t\geq 0.
$$
\end{remark}

\section{Analysis of an Auxiliary System with Viscous Regularization}\label{ecvs}
\setcounter{equation}{0}

In this section, we analyze the initial boundary value problem for the following coupled Cahn-Hilliard system with viscous regularization in chemical potentials (cf. \cite{N1988}):
\begin{align}
     & \partial_t u =\Delta \mu && \mathrm{in} \ \Omega \times (0,T),\label{eq:vs1}\\
     & \mu = \alpha \partial_t u - \Delta u + \partial_u F(u,v) && \mathrm{in} \ \Omega \times (0,T), \label{eq:vs2}\\
     & \partial_t v + \sigma (v-c)=\Delta \varphi && \mathrm{in} \ \Omega \times (0,T), \label{eq:vs3}\\
     & \varphi = \alpha \partial_t v - \Delta v + \partial_v F(u,v) && \mathrm{in} \ \Omega \times (0,T), \label{eq:vs4}\\
     & \partial_\mathbf{n} u = \partial_\mathbf{n} v  = \partial_\mathbf{n} \mu  = \partial_\mathbf{n} \varphi  = 0 && \mathrm{on} \ \partial \Omega \times (0,T), \label{eq:vsb} \\
     & (u,v)|_{t=0} = (u_0,v_0), && \mathrm{in} \ \Omega,\label{eq:vsi}
\end{align}
where $\alpha \in (0,1)$ stands for the coefficient of viscosity.

Let us summarize the main results.

\bp[Strong solutions]\label{ssve}
Let $\Omega \subset \mathbb{R}^d$ $(d\in\{2,3\})$ be a bounded domain with smooth boundary $\partial \Omega$. Suppose that $\alpha \in (0,1)$, $\delta_0\in (0, 1-|c|)$ and the assumptions $(\mathbf{H0})$, $(\mathbf{H1})$, $(\mathbf{H3})$ are satisfied.
Define
$$
W_{\delta_0}: = \left\{f \in W \ \big| \ \|f\|_{C(\overline{\Omega})}\leq 1-\delta_0\right\}.
$$
For any initial data $u_0, v_0 \in W_{\delta_0}$, problem \eqref{eq:vs1}--\eqref{eq:vsi}  admits a unique global strong solution $(u,v,\mu,\varphi)$ on $[0,+\infty)$ that satisfies the following properties:
\begin{align*}
    & u,v \in C([0,+\infty);W) \cap L^2_{\mathrm{uloc}}(0,+\infty;H^3(\Omega)), \\
    & \partial_t u, \partial_t v\in C([0,+\infty);H) \cap L^2_{\mathrm{uloc}}(0,+\infty;V) \cap H^1_{\mathrm{uloc}}(0,+\infty;V^*),\\
    & \mu, \varphi \in C([0,+\infty);W) \cap L^2_{\mathrm{uloc}}(0,+\infty;H^3(\Omega)) \cap H^1_{\mathrm{uloc}}(0,+\infty; V),
\end{align*}
$(u,v,\mu,\varphi)$
satisfies the equations \eqref{eq:vs1}--\eqref{eq:vs4} almost everywhere in $\Omega \times (0,+\infty)$ and $(u(0),v(0)) = (u_0,v_0)$ in $\Omega$. Moreover, there exists $\delta_1\in (0,\delta_0]$ such that
\begin{equation}
   \|u(t)\|_{C(\overline{\Omega})}\le 1 - \delta_1,\quad
    \|v(t)\|_{C(\overline{\Omega})} \le 1 - \delta_1, \quad \forall\, t \geq 0.
    \label{sep:rstr}
\end{equation}
The constant $\delta_1$ depends on
 $\Psi(u_0,v_0)$, $\|u_0\|_{H^2(\Omega)}$, $\|v_0\|_{H^2(\Omega)}$, $\| \nabla \mu (0) \|$, $\|\nabla \varphi (0) \|$, $ \|\partial_t u(0) \|$, $\| \partial_t v(0) \|$, $\Omega$, $\overline{u_0}$, $\overline{v_0}$, $\delta_0$, and coefficients of the system.
\ep

\begin{remark}\rm
In view of the equations \eqref{eq:vs1}--\eqref{eq:vs4}, for any $\alpha\in (0,1)$, the initial values of $\mu$, $\varphi$, $\partial_t u$, $\partial_t v$ can be determined by
\begin{align*}
    & \mu(0) = \left(\mathrm{I}-\alpha \Delta\right)^{-1} \left[-\Delta u_0 + \partial_uF(u_0,v_0)\right] \in W, \\
    & \varphi(0) = \left(\mathrm{I}-\alpha \Delta\right)^{-1} \left[-\Delta v_0 - \alpha\sigma(v_0-c) + \partial_vF(u_0,v_0)\right]\in W, \\
    & \partial_t u (0) = \Delta \mu(0) \in H, \\
    & \partial_t v (0) = \Delta \varphi(0) - \sigma (v_0-c) \in H.
\end{align*}
\end{remark}


\bp[Weak solutions] \label{wsve}
Let $\Omega \subset \mathbb{R}^d$ $(d\in\{2,3\})$ be a bounded domain with smooth boundary $\partial \Omega$. Suppose that $\alpha \in (0,1)$ and the assumptions $(\mathbf{H0})$, $(\mathbf{H1})$, $(\mathbf{H3})$ are satisfied.

(1) \emph{Existence and uniqueness}. For any initial data $u_0, v_0 \in V$ such that $F(u_0,v_0) \in L^1(\Omega)$ and $\overline{u_0},\overline{v_0} \in (-1,1)$,
problem \eqref{eq:vs1}--\eqref{eq:vsi} admits a unique global weak solution  $(u,v,\mu,\varphi)$ on $[0,+\infty)$ that satisfies the following properties:
\begin{align*}
& u, v \in C([0,+\infty);V)
\cap L^2_{\mathrm{uloc}}(0,+\infty;W) \cap H^1_{\mathrm{uloc}}(0,+\infty;V^*),\\
& \sqrt{\alpha}\partial_t u,\ \sqrt{\alpha}\partial_t v \in L^2_{\mathrm{uloc}}(0,+\infty;H),\\
& u, v \in L^\infty(\Omega\times (0,+\infty))\ \ \text{with }\ |u(x,t)|, |v(x,t)|<1\ \ \text{a.e. in}\ \Omega\times (0,+\infty),\\
&\mu, \varphi \in  L^2_{\mathrm{uloc}}(0,+\infty;V),\quad \widehat{S}'_{(u)}(u),\ \widehat{S}'_{(v)}(v)  \in L^2_{\mathrm{uloc}}(0,+\infty;H),
\end{align*}
and
\begin{align*}
    & \left\langle \partial_t u, \eta \right\rangle + (\nabla \mu, \nabla\eta) = 0,&&\forall\,\eta\in V, \ \ \text{a.e. in}\ (0,+\infty),  \\
    & \left\langle \partial_t v, \eta \right\rangle + \sigma(v-c,\eta) + (\nabla \varphi, \nabla \eta) = 0,&&\forall\,\eta\in V, \ \ \text{a.e. in}\ (0,+\infty), \\
    &\mu= \alpha \partial_t u -\Delta u+ \partial_u F(u,v), &&
     \text{a.e. in}\ \Omega\times (0,+\infty),  \\
    &\varphi= \alpha \partial_t v-\Delta v+ \partial_v F(u,v),&&
     \text{a.e. in}\ \Omega\times (0,+\infty), 
\end{align*}
moreover, it holds
$\partial_\mathbf{n} u = \partial_\mathbf{n} v= 0$ almost everywhere on $\partial \Omega\times (0,+\infty)$, and $(u(0),v(0)) = (u_0,v_0)$ almost everywhere in $\Omega$.

(2) \emph{Instantaneous regularity}. For any $\kappa \in(0,1]$, the global weak solution obtained in (1) satisfies
\begin{align*}
& u,v \in L^{\infty}(\kappa,+\infty;W)\cap H^1_{\mathrm{uloc}}(\kappa,+\infty;V), \\
& \p_t u, \p_t v \in L^\infty(\kappa,+\infty,V^*) ,\quad \sqrt{\alpha}\partial_t u,\ \sqrt{\alpha}\partial_t v \in L^{\infty}(\kappa,+\infty;H), \\
& \mu,\varphi \in L^{\infty}(\kappa,+\infty;V) , \quad \widehat{S}'_{(u)}(u),\ \widehat{S}'_{(v)}(v) \in L^{\infty}(\kappa,+\infty;H),
\end{align*}
Furthermore, there exists a constant $\delta_2 \in (0,1)$ depending on $\Psi(u_0,v_0)$, $\Omega$, $\overline{u_0}$, $\overline{v_0}$, coefficients of the system and $\kappa$, such that
\begin{equation}\label{rsep-w}
   \|u(t)\|_{C(\overline{\Omega})}\leq 1 - \delta_2,\quad  \|v(t)\|_{C(\overline{\Omega})}\leq 1 - \delta_2, \quad \forall\, t \geq \kappa.
\end{equation}
\ep

\subsection{\textit{A priori} estimates}
\label{apro-es}
Following \cite{M2019,MZ2004}, we first provide a formal derivation of some \textit{a priori} estimates for solutions to the regularized problem \eqref{eq:vs1}--\eqref{eq:vsi}.
To this end, we make the following assumptions:
\begin{itemize}
\item For some given $m \in (0,1-|c|)$, we assume $\overline{u_0}, \overline{v_0} \in [-1+m,1-m]$.
\item We assume that the initial data are sufficiently regular, i.e., $u_0, v_0\in W$ with  $$
F(u_0,v_0)\in L^1(\Omega),\quad \widehat{S}'_{(u)}(u_0),\ \widehat{S}'_{(v)}(v_0)
\in L^2(\Omega).
$$
This implies that the initial energy $\Psi(u_0,v_0)$ is bounded, moreover, $|u_0(x)|< 1$, $|v_0(x)|< 1$ almost everywhere in $\Omega$ thanks to $\mathbf{(H1)}$.
\item We assume that the solution $(u,v,\mu,\varphi)$ is sufficiently regular and fulfills
\begin{equation}
    \|u(t)\|_{L^\infty(\Omega)}<1,\quad
    \|v(t)\|_{L^\infty(\Omega)}<1,\quad \ \forall\, t\geq 0.
    \label{rLinf}
\end{equation}
\end{itemize}

\bl[Mass relations]
Let the above assumptions be satisfied. The solution $(u,v)$ to problem \eqref{eq:vs1}--\eqref{eq:vsi} satisfies
\begin{equation}\label{eq:valuesofaveuv-r}
\overline{u}(t) = \overline{u_0}\quad \text{and}\quad
\overline{v}(t) = \overline{v_0} e^{-\sigma t}+ c\left(1-e^{-\sigma t}\right),
\quad \forall\, t\geq 0.
\end{equation}
\el
\begin{proof}
Integrating \eqref{eq:vs1}, \eqref{eq:vs3} over $\Omega$, using integration by parts and the boundary conditions for $\mu$, $\varphi$, we easily arrive at the conclusion.
\end{proof}

\bl[Dissipative estimate] \label{dev}
Let the above assumptions be satisfied.
Then we have
\begin{align}
        &\widehat{\Psi}(u(t),v(t)) \leq  \widehat{\Psi}(u_0,v_0) e^{-(1+\sigma) t} + C(1 + \sigma^2), \quad \forall\, t\geq 0, \label{eq:dissipativeestimatev-1}
        \\
        & \frac{1}{2} \int_t^{t+1} \left( \| \nabla \mu (\tau ) \|^2 + \| \nabla \varphi (\tau) \|^2 \right) \mathrm{d}\tau  + \frac{\alpha}{2} \int_t^{t+1} \left( \| \partial_t u(\tau)\|^2 + \| \partial_t v(\tau)\|^2  \right) \mathrm{d}\tau \notag\\
        &\quad \leq  \widehat{\Psi}(u_0,v_0) e^{-(1+\sigma) t} + C(1 + \sigma^2), \quad \forall\, t\geq 0, \label{eq:dissipativeestimatev-2},
\end{align}
with
 \begin{align}
        \widehat{\Psi}(u,v) & = \Psi(u,v) +\frac{1}{2} \left(\| u-\overline{u}\|_{V_0^*}^2 +
        \| v-\overline{v} \|_{V_0^*}^2 \right)
        + \frac{\alpha (1+\sigma)}{2} \left(\| u-\overline{u} \|^2
        + \| v-\overline{v}\|^2\right).
        \label{newPsi1}
    \end{align}
 Here, $C$ is a positive constant depending on $\Omega$, $m$ and the parameters of the system, however, it does not depend on $\alpha$ and $\sigma$.
\el

\begin{proof}
Like in \cite[Section 4.3]{DG2022}, testing the equations \eqref{eq:vs1}, \eqref{eq:vs2}, \eqref{eq:vs3}, \eqref{eq:vs4} by $\mu$, $\partial_t u$, $\varphi$, $\partial_t v$, respectively, adding the resultants together, we obtain (cf. also \cite[Lemma 4.1]{RS2004} for computations involving singular terms)
\begin{equation}\label{eq:dissipativedifferentiale}
    \frac{\mathrm{d}}{\mathrm{dt}} \Psi(u,v)  + \| \nabla \varphi \|^2 +\| \nabla \mu \|^2 + \alpha \left(\| \partial_t u\|^2 + \| \partial_t v\|^2\right) + \sigma \int_\Omega (v-c)\varphi\, \mathrm{d}x = 0.
\end{equation}
It follows from \eqref{eq:vs4} and \eqref{eq:valuesofaveuv-r} that
\begin{align*}
    \int_\Omega (v-c)\varphi\, \mathrm{d}x
    & = \alpha \int_\Omega (v-c)\partial_t v \,\mathrm{d}x+\frac12\|\nabla v\|^2+ \int_\Omega (v-c)\partial_v F(u,v)\,\mathrm{d}x\\
    & = \frac{\alpha}{2} \frac{\mathrm{d}}{\mathrm{d}t} \|v-\overline{v} \|^2
        + \alpha \int_\Omega  \partial_t v (\overline{v}-c) \mathrm{d}x+ \frac{1}{2} \|\nabla v \|^2
    \\
    &\quad  + \int_\Omega (v-c)\widehat{S}'_{(v)}(v)\, \mathrm{d}x + \int_\Omega  \widehat{S}_{(v)}(c)\,\mathrm{d}x + \int_{\Omega} D_1(u,v)\, \mathrm{d}x,
\end{align*}
where
\begin{equation*}
    D_1(u,v) = - \widehat{S}_{(v)}(c) +  (v-c) \big[-\theta_{0,v}v + \partial_v W(u,v)\big].
\end{equation*}
Thanks to the assumption $\mathbf{(H1)}$, that is, the convexity of $\widehat{S}_{(v)}$,  we get
\begin{equation*}
    \int_\Omega (v-c)\widehat{S}'_{(v)}(v)\, \mathrm{d}x + \int_\Omega  \widehat{S}_{(v)}(c)\,\mathrm{d}x  \ge \int_{\Omega} \widehat{S}_{(v)}(v)\, \mathrm{d}x.
\end{equation*}
As a consequence,
\begin{align}
    \int_\Omega (v-c)\varphi\, \mathrm{d}x
    & \ge  \frac{\alpha}{2} \frac{\mathrm{d}}{\mathrm{d}t} \|v-\overline{v} \|^2
    + \frac{1}{2} \|\nabla v \|^2 + \int_{\Omega} \widehat{S}_{(v)}(v)\, \mathrm{d}x \notag \\
    &\quad + \int_{\Omega} D_1(u,v)\, \mathrm{d}x  + \alpha \int_\Omega \partial_t v (\overline{v}-c)\,\mathrm{d}x.
    \label{lowbd-v1}
\end{align}
By a similar argument, we have
\begin{align}
& \int_\Omega (u-\overline{u})\mu\, \mathrm{d}x
\ge \frac{\alpha}{2} \frac{\mathrm{d}}{\mathrm{d}t} \|u-\overline{u} \|^2
 + \frac{1}{2} \|\nabla u \|^2 + \int_{\Omega} \widehat{S}_{(u)}(u)\, \mathrm{d}x + \int_{\Omega} D_2(u,v)\, \mathrm{d}x,
 \label{lowbd-u1}
\end{align}
where
$$
D_2(u,v)=  - \widehat{S}_{(u)}(\overline{u}) +  (u-\overline{u}) \big[-\theta_{0,u}u + \partial_uW(u,v)\big].
$$
Combining the above estimates, we infer from \eqref{eq:dissipativedifferentiale}  that
\begin{align}
        & \frac{\mathrm{d}}{\mathrm{d}t} \left[\Psi(u,v) + \frac{\alpha \sigma}{2} \left( \|u-\overline{u}\|^2 +\|v-\overline{v}\|^2 \right)\right]\notag \\
        &\qquad + \| \nabla \mu \|^2 + \| \nabla \varphi\|^2  + \alpha \left(\| \partial_t u\|^2 +  \| \partial_t v\|^2\right) + \sigma \Psi(u,v) \notag \\
        &\quad \le \sigma \int_{\Omega} D_3\left(u,v\right) \mathrm{d}x
        + \sigma \int_\Omega (u-\overline{u})\mu\,\mathrm{d}x
        - \alpha \sigma \int_\Omega \partial_t v (\overline{v}-c)\,\mathrm{d}x,
        \label{eq:dissipativedifferentialieb}
\end{align}
where
\begin{align*}
D_3(u,v)= -D_1(u,v)-D_2(u,v)-\frac{\theta_{0,u}}{2} u^2 -\frac{\theta_{0,v}}{2} v^2 +W(u,v).
\end{align*}
The assumption $\mathbf{(H3)}$ yields $D_3\in C^1(\mathbb{R}^2)$. Hence, it follows from \eqref{rLinf} that
$$\left|\int_{\Omega} D_3(u,v) \,\mathrm{d}x\right|\leq C,$$
where the positive constant $C$ only depends on the parameters of the problem (except $\alpha, \sigma$). Using the Poincar\'{e}-Wirtinger inequality and Young's  inequality, we deduce from \eqref{rLinf} that (see \cite{DG2022})
\begin{equation*}
    \left|\int_\Omega (u-\overline{u})\mu\,\mathrm{d}x\right|
    = \left|\int_\Omega (u-\overline{u})(\mu-\overline{\mu})\,\mathrm{d}x\right|
    \le 2\left\|\mu - \overline{\mu}\right\|_{L^1\left(\Omega\right)}
    \le   \frac{1}{2\sigma} \left\|\nabla \mu \right\|^2+ C\sigma,
\end{equation*}
where $C>0$ only depends on $\Omega$. Besides, \eqref{eq:valuesofaveuv-r} yields
\begin{align*}
& \left|\int_\Omega \partial_t v (\overline{v}-c)\,\mathrm{d}x\right|
\leq \sigma |\Omega|| \overline{v}-c|^2 \leq 4 \sigma |\Omega|.
\end{align*}
Collecting the above estimates, we infer from \eqref{eq:dissipativedifferentialieb} and the fact $\alpha \in (0,1)$ that
\begin{align}
        & \frac{\mathrm{d}}{\mathrm{d}t} \left[\Psi(u,v) + \frac{\alpha \sigma}{2} \left( \|u-\overline{u}\|^2 +\|v-\overline{v}\|^2 \right)\right]\notag \\
        &\qquad + \frac12 \| \nabla \mu \|^2 + \| \nabla \varphi\|^2  + \alpha \left(\| \partial_t u\|^2 +  \| \partial_t v\|^2\right) + \sigma \Psi(u,v) \notag \\
        &\quad \le C\sigma(1+\sigma),
        \label{eq:dissipativedifferentialie}
\end{align}
where $C>0$ is independent of $\alpha$.
Next, testing \eqref{eq:vs1}, \eqref{eq:vs3} by $\mathcal{N}(u-\overline{u})$, $\mathcal{N}(v-\overline{v})$, respectively, adding the resultants together,
we have
\begin{align}
    &\frac{1}{2}\frac{\mathrm{d}}{\mathrm{d}t}\left( \| u-\overline{u}\|_{V_0^*}^2
     + \| v-\overline{v} \|_{V_0^*}^2\right)
     +\int_\Omega (u-\overline{u})\mu\,\mathrm{d}x + \int_\Omega (v-\overline{v})\varphi\,\mathrm{d}x + \sigma \left\| v-\overline{v} \right\|_{V_0^*}^2 = 0, \notag
\end{align}
Then it follows from \eqref{lowbd-u1} and similar arguments for \eqref{lowbd-v1} that
\begin{align}
        &\frac{1}{2}\frac{\mathrm{d}}{\mathrm{d}t}\left( \| u-\overline{u} \|_{V_0^*}^2
        + \| v-\overline{v} \|_{V_0^*}^2 + \alpha\| u-\overline{u} \|^2 + \alpha \| v-\overline{v} \|^2 \right) \notag \\
        & \quad +\sigma\left\| v-\overline{v} \right\|_{V_0^*}^2 + \Psi(u,v) \leq C.
        \label{v0*estimate}
\end{align}
Combining the inequalities \eqref{eq:dissipativedifferentialie}, \eqref{v0*estimate} and using  \eqref{rLinf}, we find
    \begin{align}
        &\frac{\mathrm{d}}{\mathrm{dt}}\widehat{\Psi}(u,v) + \left(1 +\sigma \right)\widehat{\Psi}(u,v)  + \frac12 \left(\| \nabla \mu \|^2 + \| \nabla \varphi\|^2\right)  + \alpha \left(\| \partial_t u\|^2 +  \| \partial_t v\|^2\right)\notag   \\
         &\quad \le C(1+\sigma^2),
         \label{es-Gron1}
    \end{align}
where $\widehat{\Psi}(u,v)$ is defined as in \eqref{newPsi1}.
It follows from \eqref{es-Gron1} and Gronwall's lemma that
\begin{equation*}
    \widehat{\Psi}(u(t),v(t))
     \le e^{-(1+\sigma) t} \widehat{\Psi}(u_0,v_0) + C(1+\sigma^2),
     \quad \forall\, t\geq 0.
\end{equation*}
This gives \eqref{eq:dissipativeestimatev-1} with $C$ being a positive constant that only depends on $\Omega$, $m$ and parameters of the system except $\alpha$ and $\sigma$. On the other hand, from $\mathbf{(H1)}$, $\mathbf{(H3)}$ and \eqref{rLinf}, we observe that
\begin{align}
\widehat{\Psi}(u(t),v(t)) \geq -\frac12|\Omega|\left(\theta_{0,u}+\theta_{0,v}\right)- |\Omega|\max_{(r,s)\in [-1,1]^2}|W(r,s)|.
\label{low-bd1}
\end{align}
Thus, integrating \eqref{es-Gron1} on $[t,t+1]$, we can conclude the dissipative estimate
\eqref{eq:dissipativeestimatev-2}.
\end{proof}

\bl[Lower-order estimates]\label{rlow}
Let the above assumptions be satisfied.
We have
\begin{align*}
    & u,v \in   L^{\infty}(0,+\infty;V)\cap L^2_{\mathrm{uloc}}(0,+\infty;W) \cap H^1_{\mathrm{uloc}}(0,+\infty;V^*),\\
    & \sqrt{\alpha}\partial_t u,\ \sqrt{\alpha}\partial_t v \in L^2_{\mathrm{uloc}}(0,+\infty;H),\\
    & \mu,\ \varphi \in L^2_{\mathrm{uloc}}(0,+\infty;V),\\
    & \widehat{S}'_{(u)}(u),\ \widehat{S}'_{(v)}(v)  \in L^2_{\mathrm{uloc}}(0,+\infty;H),
\end{align*}
with uniform bounds with respect to $\alpha\in (0,1)$ in the corresponding spaces.
\el
\begin{proof}
Thanks to Lemma \ref{dev} and \eqref{rLinf}, we easily find that
\begin{equation}
\left\{\
\begin{aligned}
    & u,v \in   L^{\infty}(0,+\infty;V),\quad \sqrt{\alpha}\partial_t u,\ \sqrt{\alpha}\partial_t v \in L^2_{\mathrm{uloc}}(0,+\infty;H),\\
    & \nabla \mu,\ \nabla \varphi \in L^2_{\mathrm{uloc}}(0,+\infty;\bm{L}^2(\Omega)),
    \end{aligned}
    \right.
    \label{rlowes-1}
\end{equation}
are uniformly bounded with respect to $\alpha\in (0,1)$ in the corresponding spaces. By comparison in \eqref{eq:vs1} and \eqref{eq:vs3}, we also get
$$\partial_tu, \partial_t v\in H^1_{\mathrm{uloc}}(0,+\infty;V^*).$$
According to \cite{GGG2017,MZ2004}, we have the following inequality for
$\widehat{S}'_{(u)}$:
\begin{equation*}
     \left\|\widehat{S}'_{(u)}(u) \right\|_{L^1(\Omega)} \le \mathcal{R}_{(u)}(|\overline{u}|)\left[1+\int_\Omega  \left(\widehat{S}'_{(u)}(u)-\overline{\widehat{S}'_{(u)}(u)}\right) \left(u-\overline{u}\right) \,\mathrm{d}x\right],
\end{equation*}
where $\mathcal{R}_{(u)}(\cdot)$ is an increasing function.  By comparison in \eqref{eq:vs2}, we find
\begin{equation*}
    \begin{aligned}
        \left\|\widehat{S}'_{(u)}(u )-\overline{\widehat{S}'_{(u)}(u)}\right\|_{V^*}
        & \le \|\mu -\overline{\mu}\|_{V^*} + \alpha \| \partial_t u \|_{V^*} + \|\Delta u \|_{V^*} + \|D_{(u)}(u,v) -\overline{D_{(u)}(u,v)}\|_{V^*} \\
        & \le C \left(\|\nabla \mu\| + \|\nabla u\| + \| D_{(u)}(u,v)\|\right)\\
        & \le C \left(1 + \|\nabla \mu\|\right),
    \end{aligned}
\end{equation*}
where
\begin{equation*}
    D_{(u)}(u,v) =\partial_u W(u,v)-\theta_{0,u} u.
\end{equation*}
Therefore, it holds
\begin{equation*}
        \left\|\widehat{S}'_{(u)}(u) \right\|_{L^1(\Omega)} \le C \left(1  +\|\nabla \mu\|\right),
\end{equation*}
where $C >0$ depends on $\Psi(u_0,v_0)$, $\Omega$, $m$, and coefficients of the system except $\alpha$.
In a similar manner, we obtain the following estimate for $\widehat{S}'_{(v)}(v)$:
\begin{equation*}
        \left\|\widehat{S}'_{(v)}(v) \right\|_{L^1(\Omega)} \le C \left(1 +\|\nabla \varphi\|\right).
\end{equation*}
The above two estimates yield
\begin{align}
    &|\overline{\mu}| = \left|\overline{\partial_u F(u,v)}\right| \le C \left(1   +\|\nabla \mu\|\right),
    \label{eq:avemu}\\
    &|\overline{\varphi}| = \left|\overline{\partial_vF(u,v)} +\alpha\overline{\partial_t v}\right| \le C\left(1 + \|\nabla \varphi\|\right).
    \label{eq:avepsi}
\end{align}
Hence, by the Poincar\'e-Wirtinger inequality, we find
\begin{align}
\mu,\ \varphi \in L^2_{\mathrm{uloc}}(0,+\infty;V).
\label{rlowes-2}
\end{align}
Let us rewrite the equations \eqref{eq:vs2}, \eqref{eq:vs4} as
\begin{align*}
& -\Delta u + \widehat{S}'_{(u)}(u) =\mu - \alpha \partial_t u + \theta_{0,u}u - \partial_u W(u,v),\\
& -\Delta v + \widehat{S}'_{(v)}(u) =\varphi - \alpha \partial_t v + \theta_{0,v}v - \partial_v W(u,v).
\end{align*}
Keeping \eqref{rlowes-1}, \eqref{rlowes-2} in mind, we can apply Lemma \ref{nonlinearelliptic} to conclude that
$$
u,v \in L^2_{\mathrm{uloc}}(0,+\infty;W),\quad \widehat{S}'_{(u)}(u),\ \widehat{S}'_{(v)}(v)  \in L^2_{\mathrm{uloc}}(0,+\infty;H),
$$
with uniform bounds with respect to $\alpha\in (0,1)$.
\end{proof}

\bl[Higher-order estimates]\label{hoev}
Let the above assumptions be satisfied. For any given $\kappa \in (0,1]$, we have
\begin{equation}\label{eq:hdevs}
    \begin{aligned}
        & \| \nabla \mu (t)\|^2 + \| \nabla \varphi (t) \|^2 + \alpha \left(\| \partial_t u(t)\|^2 + \| \partial_t v(t)\|^2\right)
          + \int_t^{t+1} \left(\| \nabla \partial_t u(\tau) \|^2 + \| \nabla \partial_t v(\tau) \|^2\right) \mathrm{d\tau} \\
        &\quad \le \chi_{[0,1]}\left(1 - \frac{t}{\kappa}\right) \mathcal{Q}_1 + \mathcal{Q}_2,\quad \forall\, t\geq 0,
    \end{aligned}
\end{equation}
where $\chi_{[0,1]}$ denotes the indicator function of $[0,1]$, the positive constant $\mathcal{Q}_1$ depends on $\Psi(u_0,v_0)$, $\|u_0\|_{H^2(\Omega)}$, $\|v_0\|_{H^2(\Omega)}$, $\| \nabla \mu (0) \|$, $\|\nabla \varphi (0) \|$, $ \|\partial_t u(0) \|$, $\| \partial_t v(0) \|$, $\Omega$, $m$ and coefficients of the system, while the positive constant $\mathcal{Q}_2$ depends on $\Psi(u_0,v_0)$, $\Omega$, $m$, $\kappa$ and coefficients of the system (except $\alpha$) . As a consequence, it holds
\begin{align*}
& \partial_t u,\ \partial_t v \in L^2_{\mathrm{uloc}}(\kappa,+\infty;V)\cap L^\infty(\kappa,+\infty,V^*),\\
& \sqrt{\alpha}\partial_t u,\ \sqrt{\alpha}\partial_t v \in L^{\infty}(\kappa,+\infty;H),\qquad
 \mu,\varphi \in L^{\infty}(\kappa,+\infty;V),\\
 &   u,v \in L^{\infty}(\kappa,+\infty;W),\qquad
\widehat{S}'_{(u)}(u),\ \widehat{S}'_{(v)}(u) \in L^{\infty}(\kappa,+\infty;H),
\end{align*}
with uniform bounds with respect to $\alpha\in (0,1)$ in the corresponding spaces.
\el
\begin{proof}
Like in \cite{GGM2017}, testing \eqref{eq:vs1}, \eqref{eq:vs3} with $\partial_t \varphi$, $\partial_t\varphi$, respectively, we obtain
\begin{align}
    & \frac{1}{2} \frac{\mathrm{d}}{\mathrm{dt}} \|\nabla \mu \|^2
    + \int_\Omega \partial_t u\partial_t\mu\,\mathrm{d}x = 0 ,\label{muttestu}\\
    & \frac{1}{2} \frac{\mathrm{d}}{\mathrm{dt}} \|\nabla \varphi \|^2
    + \int_\Omega \partial_t v\partial_t\varphi\,\mathrm{d}x
    + \sigma\int_\Omega (v-c)\partial_t\varphi\,\mathrm{d}x = 0.\label{phittestv}
\end{align}
Using \eqref{eq:vs2}, \eqref{rLinf}, $\mathbf{(H1)}$, $\mathbf{(H3)}$ and integration by parts, we find
    \begin{align}
        \int_\Omega \partial_t u\partial_t\mu\,\mathrm{d}x
        & = \frac{\alpha}{2} \frac{\mathrm{d}}{\mathrm{d}t} \|\partial_t u \|^2  + \|\nabla\partial_t u  \|^2
        + \int_\Omega \frac{\partial^2 F}{\partial u^2} (\partial_t u)^2\,\mathrm{d}x
        + \int_\Omega \frac{\partial^2 F}{ \partial u \partial v}  \partial_t v \partial_t u  \,\mathrm{d}x \notag \\
        & \ge \frac{\alpha}{2} \frac{\mathrm{d}}{\mathrm{d}t} \|\partial_t u \|^2  + \|\nabla\partial_t u  \|^2  - C\left(\|\partial_t u\|^2+ \|\partial_t v\|^2\right),
        \label{eq:utmut}
    \end{align}
and in a similar manner,
    \begin{align}
        \int_\Omega \partial_t v\partial_t\varphi\,\mathrm{d}x
        & = \frac{\alpha}{2} \frac{\mathrm{d}}{\mathrm{d}t} \|\partial_t v \|^2  + \|\nabla\partial_t v  \|^2
       +  \int_\Omega \frac{\partial^2 F}{ \partial v \partial u}  \partial_t u \partial_t v  \,\mathrm{d}x
        + \int_\Omega \frac{\partial^2 F}{\partial v^2} (\partial_t v)^2\,\mathrm{d}x
          \notag \\
        & \ge \frac{\alpha}{2} \frac{\mathrm{d}}{\mathrm{d}t} \|\partial_t v \|^2  + \|\nabla\partial_t v  \|^2  - C\left(\|\partial_t u\|^2+ \|\partial_t v\|^2\right).
        \label{eq:vtpsit}
    \end{align}
Besides, it follows from \eqref{eq:vs4} that
\begin{align}
        \int_\Omega (v-c)\partial_t\varphi\,\mathrm{d}x
        & = \frac{\mathrm{d}}{\mathrm{d}t} \int_\Omega(v-c)\varphi \,\mathrm{d}x
        - \int_\Omega \varphi \partial_t v  \,\mathrm{d}x \notag\\
        & = \frac{\mathrm{d}}{\mathrm{d}t} \left(\int_\Omega (v-c)\varphi \,\mathrm{d}x
        - \frac{1}{2}\|\nabla v \|^2 - \int_{\Omega} \widehat{S}_{(v)}(v) \,\mathrm{d}x\right) \notag\\
        &\quad -\int_\Omega \left(\partial_v W(u,v)-\theta_{0,v}v\right) \partial_t v\,\mathrm{d}x
         - \alpha \|\partial_t v\|^2\notag\\
         &\geq \frac{\mathrm{d}}{\mathrm{d}t} \left(\int_\Omega (v-c)\varphi \,\mathrm{d}x
        - \frac{1}{2}\|\nabla v \|^2 - \int_{\Omega} \widehat{S}_{(v)}(v) \,\mathrm{d}x\right) - 2\|\partial_t v\|^2-C,
         \label{eq:v-cpsit}
    \end{align}
where in the last step we use Young's inequality and the fact $\alpha \in (0,1)$.
Substituting \eqref{eq:utmut}, \eqref{eq:vtpsit}  and \eqref{eq:v-cpsit} into the sum of \eqref{muttestu} and \eqref{phittestv}, by interpolation and Young's inequality, we obtain
    \begin{align}
        & \frac{\mathrm{d}}{\mathrm{d}t}\mathcal{H}(t)  + \|\nabla \partial_t u\|^2 + \|\nabla \partial_t v\|^2 \notag  \\
        &\quad \le  C\left(\|\partial_t u\|^2+ \|\partial_t v\|^2\right)+C\notag  \\
        &\quad \leq \frac12 \|\nabla \partial_t u\|^2 + C \|\partial_t u\|_{V^*}^2
        + \frac12 \|\nabla \partial_t v\|^2 + C \|\partial_t (v-\overline{v})\|_{V^*}^2+ C\|\partial_t\overline{v}\|^2+C\notag\\
        &\quad \leq \frac12 \left(\|\nabla \partial_t u\|^2 + \|\nabla \partial_t v\|^2\right) + C \left(\|\partial_t u\|_{V^*}^2
        + \|\partial_t v\|_{V^*}^2\right)+ C,
        \label{eq:highorderdissipativedifferentialie}
    \end{align}
where
\begin{align}
\mathcal{H}(t)& = \frac{1}{2}
        \left( \| \nabla \mu \|^2 + \| \nabla \varphi\|^2 + \alpha \| \partial_t u\|^2 + \alpha \| \partial_t v\|^2\right) \notag\\
         & \quad + \sigma   \left(\int_\Omega (v-c)\varphi \,\mathrm{d}x
        - \frac{1}{2}\|\nabla v \|^2 - \int_{\Omega} \widehat{S}_{(v)}(v) \,\mathrm{d}x\right),
        \label{H1}
\end{align}
and the constant $C>0$ is independent of $\alpha$.

Recalling \eqref{lowbd-v1}, we infer from \eqref{rLinf}, \eqref{eq:valuesofaveuv-r} that
    \begin{align}
        & \int_\Omega (v-c)\varphi\,\mathrm{d}x - \frac{1}{2}\|\nabla v\|^2 - \int_{\Omega} \widehat{S}_{(v)}(v) \,\mathrm{d}x\notag \\
         &\quad \ge \frac{\alpha}{2} \frac{\mathrm{d}}{\mathrm{d}t} \|v-\overline{v} \|^2
         + \int_{\Omega} D_1(u,v)\, \mathrm{d}x  + \alpha \int_\Omega \partial_t v (\overline{v}-c)\,\mathrm{d}x\notag \\
         &\quad \ge -\frac{\alpha}{4\sigma} \|\partial_t v\|^2 - C,
 \label{eq:v-cpsi}
    \end{align}
where $C>0$ is independent of $\alpha\in(0,1)$. On the other hand, from \eqref{rLinf}, \eqref{eq:avepsi} and the Poincar\'e-Wirtinger inequality, we find the upper bound
\begin{align}
 & \int_\Omega (v-c)\varphi\,\mathrm{d}x - \frac{1}{2}\|\nabla v\|^2 - \int_{\Omega} \widehat{S}_{(v)}(v) \,\mathrm{d}x
 \leq  \|v-c\|\|\varphi\|
 \leq  C\left(1 + \|\nabla \varphi\|\right).
         \label{eq:v-cpsib}
\end{align}
According to \eqref{H1}, \eqref{eq:v-cpsi}, \eqref{eq:v-cpsib}, there exist some constants $C_1,C_2>0$ independent of $\alpha$ such that
\begin{align}
& \frac{1}{4}
        \left( \| \nabla \mu \|^2 + \| \nabla \varphi\|^2 + \alpha \| \partial_t u\|^2 + \alpha \| \partial_t v\|^2\right)\notag  \\
&\quad \leq \mathcal{H}(t)+C_1 \leq   \| \nabla \mu \|^2 + \| \nabla \varphi\|^2 + \alpha \| \partial_t u\|^2 + \alpha \| \partial_t v\|^2+C_2.
\label{l-h-bdH1}
\end{align}

Now for any given $t \in (0,2]$, integrating \eqref{eq:highorderdissipativedifferentialie} over $[0,t]$, we get
    \begin{align}
        & \mathcal{H}(t) + \frac12\int_0^t \left(\|\nabla \partial_t u(\tau)\|^2 + \|\nabla \partial_t v(\tau)\|^2\right)\,\mathrm{d}\tau \notag\\
        &\quad \le \mathcal{H}(0)  + C \int_0^t \left(\|\partial_t u(\tau)\|_{V^*}^2
        + \|\partial_t v(\tau)\|_{V^*}^2\right)\,\mathrm{d}\tau + C.
        \label{eq:inthighorderdissipativedifferentialie}
    \end{align}
Using the lower-order estimates obtained in Lemma \ref{rlow}, we infer from  \eqref{l-h-bdH1} and \eqref{eq:inthighorderdissipativedifferentialie}  that
\begin{equation*}
    \begin{aligned}
         &\frac{1}{4}
        \left( \| \nabla \mu(t) \|^2 + \| \nabla \varphi(t)\|^2 + \alpha \| \partial_t u(t)\|^2 + \alpha \| \partial_t v(t)\|^2\right) \notag \\
        &\qquad + \frac12\int_0^t \left(\|\nabla \partial_t u(\tau)\|^2 + \|\nabla \partial_t v(\tau)\|^2\right)\,\mathrm{d}\tau
        \le C,\quad \forall\, t\in [0,2],
    \end{aligned}
\end{equation*}
which gives
    \begin{align}
        &  \| \nabla \mu (t)\|^2 + \| \nabla \varphi (t) \|^2 + \alpha \left(\| \partial_t u(t)\|^2 + \| \partial_t v(t) \|^2\right)
        +   \int_t^{t+1} \left(\| \nabla \partial_t u(\tau) \|^2 + \| \nabla \partial_t v(\tau) \|^2\right) \mathrm{d}\tau \notag \\
        &\quad \le C, \quad \forall\, t \in [0,1],
        \label{eq:hdevsfinite}
    \end{align}
where $C>0$ depends on $\Psi(u_0,v_0)$, $\|u_0\|_{H^2(\Omega)}$, $\|v_0\|_{H^2(\Omega)}$, $\| \nabla \mu (0) \|$, $\|\nabla \varphi (0) \|$, $ \|\partial_t u(0) \|$, $\| \partial_t v(0) \|$, $\Omega$, $m$, $\sigma$ and also $\alpha$.

Next, for any given $\kappa \in (0,1]$, we derive higher-order estimates on the infinite interval $[\kappa,+\infty)$. Thanks to Lemma \ref{rlow} and \eqref{l-h-bdH1}, it holds
$$
\int_t^{t+\kappa}\left(\mathcal{H}(\tau)+C_1\right)\,\mathrm{d}\tau
\leq \int_t^{t+1}\left(\mathcal{H}(\tau)+C_1\right)\,\mathrm{d}\tau
\leq C,\quad \forall\, t\geq 0,
$$
where $C>0$ depends on $\Psi(u_0,v_0)$, $\Omega$, $m$, $\sigma$, but are independent of $\alpha\in (0,1)$ and $\kappa\in (0,1]$. Thus, we can infer from the above estimate, Lemma \ref{rlow}, \eqref{eq:highorderdissipativedifferentialie} and the uniform Gronwall lemma (see Lemma 1.1, Chap. III of Ref. \cite{Temam1997}) that
\begin{align}
\mathcal{H}(t)+C_1 \leq C\left(1+\frac{1}{\kappa}\right),\quad \forall\, t\geq \kappa,\notag
\end{align}
which yields
\begin{align}
        &  \| \nabla \mu (t)\|^2 + \| \nabla \varphi (t) \|^2 + \alpha \left(\| \partial_t u(t)\|^2 + \| \partial_t v(t) \|^2\right) \le C\left(1+\frac{1}{\kappa}\right), \quad \forall\, t\geq \kappa.
        \label{eq:hdevsfinite-b}
    \end{align}
Furthermore, integrating \eqref{eq:highorderdissipativedifferentialie} in time, we get
  \begin{align}
        & \int_t^{t+1} \left(\| \nabla \partial_t u(\tau) \|^2 + \| \nabla \partial_t v(\tau) \|^2\right) \mathrm{d}\tau
         \le C\left(1+\frac{1}{\kappa}\right),\quad \forall\, t\geq \kappa.
        \label{eq:hdevsfinite-c}
  \end{align}

Combining the estimates \eqref{eq:hdevsfinite}, \eqref{eq:hdevsfinite-b} and \eqref{eq:hdevsfinite-c}, we arrive at the conclusion \eqref{eq:hdevs}.
In particular, it follows from \eqref{eq:avemu}, \eqref{eq:avepsi}, \eqref{eq:hdevsfinite-b} and the Poincar\'e-Wirtinger inequality that $\mu,\varphi \in L^{\infty}(\kappa,+\infty;V)$ for any given $\kappa\in (0,1]$.
By comparison in \eqref{eq:vs1} and \eqref{eq:vs3}, we have $\partial_t u, \partial_t v\in L^\infty(\kappa,+\infty,V^*)$. On the other hand, by Lemma \ref{nonlinearelliptic} (see \eqref{eq:ellipticestimate}), we obtain $u,v \in L^{\infty}(\kappa,+\infty;W)$ and $\widehat{S}'_{(u)},\ \widehat{S}'_{(v)} \in L^{\infty}(\kappa,+\infty;H)$ with uniform bounds with respect to $\alpha\in (0,1)$.
\end{proof}

\bl[Strict separation property]\label{strictseparationvs}
Let the above assumptions be satisfied.

(1) Assume in addition, $u_0,v_0 \in W_{\delta_0}$ for some $\delta_0\in (0,m]$. Then we have
\begin{equation}\label{rsep-1}
   \|u(t)\|_{L^{\infty}(\Omega)}\leq 1 - \delta_1,\quad  \|v(t)\|_{L^{\infty}(\Omega)}\leq 1 - \delta_1, \quad  \forall\, t \geq 0,
\end{equation}
where the constant $\delta_1\in (0,1)$ depends on
 $\Psi(u_0,v_0)$, $\|u_0\|_{H^2(\Omega)}$, $\|v_0\|_{H^2(\Omega)}$, $\| \nabla \mu (0) \|$, $\|\nabla \varphi (0) \|$, $ \|\partial_t u(0) \|$, $\| \partial_t v(0) \|$, $\Omega$, $m$, $\delta_0$, and coefficients of the system.

(2) For any given $\kappa\in (0,1]$, there exists $\delta_2\in (0,1)$ depending on $\Psi(u_0,v_0)$, $\Omega$, $m$, coefficients of the system and $\kappa$, such that
\begin{equation}\label{rsep-2}
   \|u(t)\|_{L^{\infty}(\Omega)}\leq 1 - \delta_2,\quad  \|v(t)\|_{L^{\infty}(\Omega)}\leq 1 - \delta_2, \quad \forall\, t \geq \kappa,
\end{equation}
\el
\begin{remark}\rm
We emphasize that both constants $\delta_1$, $\delta_2$ in Lemma \ref{strictseparationvs} depend on the viscous parameter $\alpha\in (0,1)$.
\end{remark}

\begin{proof}
The proof follows the idea in \cite{MZ2004}, that is, for every $\alpha \in (0,1)$ one is allowed to take advantage of the comparison principle for second-order parabolic equations. To this end, we note that \eqref{eq:vs2} can be written as
\begin{equation}\label{eq:parabolicu}
    \alpha \partial_t u - \Delta u + \widehat{S}'_{(u)}(u) = h_{(u)},
\end{equation}
with
$$h_{(u)} = \theta_{0,u} u - \partial_u W(u,v) + \mu.$$

(1) It follows from Lemma \ref{rlow}, Lemma \ref{hoev}, \eqref{rLinf}, \eqref{eq:avemu}, $\mathbf{(H1)}$, $\mathbf{(H3)}$ and the Sobolev embedding $H^2(\Omega)\hookrightarrow L^\infty(\Omega)$ that
\begin{align*}
\|h_{(u)}(t)\|_{L^\infty(\Omega)}
&\leq  \|\theta_{0,u} u\|_{L^\infty(\Omega)}
+ \|\partial_u W(u,v)\|_{L^\infty(\Omega)} + C\|\mu\|_{H^2(\Omega)}\\
&\leq C+ C (|\overline{\mu}|+\|\nabla \mu\| + \|\partial_t u\|)\\
&\leq \widetilde{C},\quad \forall\, t\geq 0,
\end{align*}
where $\widetilde{C}>0$ depends on $\Psi(u_0,v_0)$, $\|u_0\|_{H^2(\Omega)}$, $\|v_0\|_{H^2(\Omega)}$, $\| \nabla \mu (0) \|$, $\|\nabla \varphi (0) \|$, $ \|\partial_t u(0) \|$, $\| \partial_t v(0) \|$, $\Omega$, $m$, $\sigma$, $\alpha$. Consider the following auxiliary ODE systems for $y_+$ and $y_{-}$:
\begin{align*}
 \alpha \partial_t y_\pm +  \widehat{S}'_{(u)}(y_\pm) = \pm \widetilde{C},\qquad y_\pm(0)=\pm\|u_0\|_{L^\infty(\Omega)}\in (-1,1).
\end{align*}
 Since $u_0\in W_{\delta_0}$, we can apply \cite[Proposition A.3]{MZ2004} to conclude
 $$
 -1+\delta_{1,u} \leq y_{-}(t)\leq 0 \leq y_+(t)\leq 1-\delta_{1,u},\quad \forall\, t\geq 0,
 $$
 where $\delta_{1,u}\in (0,\delta_0]$ depends on $\delta_0$ and $\widetilde{C}$.
 On the other hand, by the comparison principle for \eqref{eq:parabolicu}, we get
 $$
  y_{-}(t)\leq u(x,t) \leq y_+(t),\quad \forall\, (x,t)\in \Omega\times[0,+\infty).
 $$
 The above facts yield the strict separation for $u$ on $[0, +\infty)$. By a similar argument, we can obtain the strict separation property for $v$ with distance $\delta_{1,v}\in (0,\delta_0]$. Taking $\delta_1=\min\{\delta_{1,u}, \delta_{1,v}\}$, we conclude \eqref{rsep-1}.

 (2) To prove \eqref{rsep-2}, we just note that for any given $\kappa\in (0,1]$, it holds
 \begin{align*}
\|h_{(u)}(t)\|_{L^\infty(\Omega)}
&\leq  \|\theta_{0,u} u\|_{L^\infty(\Omega)}
+ \|\partial_u W(u,v)\|_{L^\infty(\Omega)} + C\|\mu\|_{H^2(\Omega)}\\
&\leq C+ C (|\overline{\mu}|+\|\nabla \mu\| + \|\partial_t u\|)\\
&\leq \widetilde{C}_\kappa,\quad \forall\, t\geq \frac{\kappa}{2},
\end{align*}
where $\widetilde{C}_\kappa>0$ depends on $\Psi(u_0,v_0)$, $\Omega$, $m$, $\sigma$, $\alpha$ and $\kappa$. Consider the following auxiliary ODE systems for $y_+$ and $y_{-}$:
\begin{align*}
 \alpha \partial_t y_\pm +  \widehat{S}'_{(u)}(y_\pm) = \pm \widetilde{C}_\kappa,\qquad y_\pm(0)=\pm\left\|u\left(\frac{\kappa}{2}\right)\right\|_{L^\infty(\Omega)}\in(-1,1).
\end{align*}
From \eqref{rLinf} and \cite[Corollary A.1]{MZ2004}, there exists a constant $\delta_2\in (0,1)$ depending on $\kappa$ and $\widetilde{C}_\kappa$ such that
  $$
 -1+\delta_2 \leq y_{-}(t)\leq 0 \leq y_+(t)\leq 1-\delta_2,\quad \forall\, t\geq \kappa.
 $$
Using the comparison principle again, we arrive at the conclusion \eqref{rsep-2}.
\end{proof}

\subsection{Well-posedness of the auxiliary problem}

First, we present a result on the continuous dependence of weak solutions to problem \eqref{eq:vs1}--\eqref{eq:vsi} with respect to the initial data, which also implies the uniqueness of weak solutions. Since the viscous terms $\alpha\partial_t u$, $\alpha\partial_t v$ do not bring any trouble for the argument as in \cite[Section 4.4]{DG2022}, we omit the detailed proof here.

\bp[Continuous dependence]\label{ssvc}
Let $d = 2,3$. Given $R\ge 0$, $T\in (0,+\infty)$ and $m \in (|c|,1)$, there exists a positive constant $C$ depending on $m$, $R$, $T$ such that, for any weak solutions $(u_1,v_1)$, $(u_2,v_2)$ to problem \eqref{eq:vs1}--\eqref{eq:vsi} on $[0,T]$ originating from the initial data $(u_{01},v_{01})$, $(u_{02},v_{02})$ satisfying $\Psi(u_{0i},v_{0i})\le R$ and $|\overline{u_{0i}}|$, $|\overline{v_{0i}}|\le m$ $(i=1,2)$, the continuous dependence estimate
\begin{equation*}
    \begin{aligned}
        & \| u_1(t) - u_2(t)\|_{V^*}^2 + \| v_1(t) - v_2(t)\|_{V^*}^2
        +  \int_0^T \| u_1(t) - u_2(t)\|_{V}^2 \,\mathrm{d}t
        + \int_0^T \| v_1(t) - v_2(t) \|_{V}^2 \,\mathrm{d}t  \\
        & \quad \le  C \left(\| u_{01} - u_{01} \|_{V^*}^2 + \| v_{01} - v_{02}\|_{V^*}^2 + |\overline{u_{01}}-\overline{u_{02}}| + |\overline{v_{01}}-\overline{v_{02}}|\right)
    \end{aligned}
\end{equation*}
holds for every $t\in \left[0,T\right]$.
\ep

Next, we prove Proposition \ref{ssve} on the existence and uniqueness of a global strong solution to problem \eqref{eq:vs1}--\eqref{eq:vsi}.

\medskip

\noindent \textbf{Proof of Proposition \ref{ssve}.} We extend the argument in \cite{MZ2004} to the coupled system \eqref{eq:vs1}--\eqref{eq:vsi}.

\textit{Step 1. Local well-posedness.}
For any $\delta>0$, we consider the following approximate problem
\begin{equation}
\left\{
\begin{aligned}
     & \partial_t u =\Delta \mu &&\qquad\quad   \mathrm{in} \ \Omega \times (0,T), \\
     & \mu = \alpha \partial_t u - \Delta u + \partial_u F_{\delta}(u,v) &&\qquad \quad  \mathrm{in} \ \Omega \times (0,T),  \\
     & \partial_t v + \sigma (v-c)=\Delta \varphi &&\qquad\quad   \mathrm{in} \ \Omega \times (0,T), \\
     & \varphi = \alpha \partial_t v - \Delta v + \partial_v F_{\delta}(u,v) && \qquad \quad \mathrm{in} \ \Omega \times (0,T),  \\
     & \partial_\mathbf{n} u = \partial_\mathbf{n} v  = \partial_\mathbf{n} \mu  = \partial_\mathbf{n} \varphi  = 0 && \qquad \quad \mathrm{on} \ \partial \Omega \times (0,T), \\
     & (u,v)|_{t=0} = (u_0,v_0), && \qquad \quad \mathrm{in} \ \Omega.
\end{aligned}
\right.
\label{eq:as1}
\end{equation}
Here, we introduce a cut-off of the nonlinear function $F$ (see \eqref{eq:F}) such that
$$
F_{\delta}(u,v) = S_{\delta}\left(u;\theta_u,\theta_{0,u}\right) + S_{\delta}\left(v;\theta_v,\theta_{0,v}\right) + W_{\delta}(u,v),
$$
where
\begin{equation}\notag
    S_{\delta}\left(s;\theta_r,\theta_{0,r}\right) = \left\{
        \begin{aligned}
            & S\left(s;\theta_r,\theta_{0,r}\right) \xi_\delta(s), &&\ s \in (-1,1), \\
            & 0, &&\ s \in \mathbb{R}\setminus(-1,1),
        \end{aligned}
        \right.
\end{equation}
 for $r\in\{u,v\}$ and
\begin{equation}\notag 
    W_{\delta}(u,v) = W(u,v)\xi_\delta(u)\xi_\delta(v).
\end{equation}
The smooth cutoff function $\xi_\delta:\mathbb{R}\to \mathbb{R}$ is given by
$$
\xi(s) = \left(\chi_{\left[-1+\frac{\delta}{2},1-\frac{\delta}{2}\right]}(\cdot)* \eta_{\frac{\delta}{4}}(\cdot)\right)(s),\quad \forall\,s \in \mathbb{R},
$$
where $*$ means the usual convolution and $\{\eta_{\epsilon}\}_{\epsilon>0}$ is a family of mollifiers on $\mathbb{R}$. It is straightforward to verify that $F_{\delta}(\cdot,\cdot)$ and its partial derivatives of order not larger than two are uniformly bounded on any bounded set in $\mathbb{R}^2$. That is, $F_{\delta}(\cdot,\cdot) \in BUC^2(\mathbb{R}\times\mathbb{R})$. Moreover,
$$
F_{\delta}(u,v)= F(u,v),\qquad \forall\, u,v\in [-1+\delta,1-\delta].
$$

From the assumption that $u_0,v_0 \in W_{\delta_0}$, we have $\overline{u_0}, \overline{v_0}\in [-1+\delta_0,1-\delta_0]$ and thus can take $m=\delta_0$. Let $\delta_1\in(0,\delta_0]$ be the strict separation constant determined by Lemma \ref{strictseparationvs}. Then we consider problem \eqref{eq:as1} with the nonlinear term $F_{\delta_1/2}(u,v)$, subject to the same initial data $(u_0,v_0)$. It follows from Proposition \ref{ssae} (see appendix) that this auxiliary problem admits a unique global strong solution, which we denote by $(u^\sharp,v^\sharp,\mu^\sharp, \varphi^\sharp)$. Since $u^\sharp,v^\sharp\in C([0,T];W)$, from the Sobolev embedding $H^2(\Omega) \hookrightarrow C(\overline{\Omega})$ that holds in two and three dimensions, there exists some positive $T_*>0$ such that
\begin{align}
\| u^\sharp(t)\|_{C(\overline{\Omega})}\leq 1- \frac{3}{4}\delta_1,\quad
\| v^\sharp(t) \|_{C(\overline{\Omega})}\leq 1- \frac{3}{4}\delta_1, \quad  \forall \, t \in [0,T_*].
\label{sep-loc1}
\end{align}
The strict separation property \eqref{sep-loc1} together with the definition of $F_{\delta_1/2}$ yields
\begin{align*}
    & \partial_u F_{\delta_1/2}(u^\sharp (t),v^\sharp(t)) = \partial_u F(u^\sharp(t),v^\sharp(t)),\quad  \\
    & \partial_v F_{\delta_1/2}(u^\sharp(t),v^\sharp(t)) = \partial_v F(u^\sharp(t),v^\sharp(t)),
\end{align*}
for all $t\in [0,T_*]$.
As a consequence, $(u^\sharp,v^\sharp)$ is indeed a strong solution to the original problem \eqref{eq:vs1}--\eqref{eq:vsi} on $[0,T_*]$. This gives the existence of a local strong solution $(u,v)=(u^\sharp,v^\sharp)$  to problem \eqref{eq:vs1}--\eqref{eq:vsi} (with corresponding chemical potentials $(\mu, \varphi)$). According to Proposition \ref{ssvc}, this local strong solution is unique on $[0,T_*]$.
\smallskip

 \textit{Step 2. Global existence.} It is obvious that the local strong solution $(u,v)$ to problem \eqref{eq:vs1}--\eqref{eq:vsi} can be (uniquely) extended beyond the finite interval $[0,T_*]$. Next, we show that $(u,v)$ is indeed global. To this end, let $T_{\mathrm{max}}$ be the maximal existence time for $(u,v)$ such that
 $$
 \| u(t)\|_{C(\overline{\Omega})}<1,\quad
\| v(t) \|_{C(\overline{\Omega})}<1, \quad  \forall \, t \in [0,T_{\mathrm{max}}).
 $$
 Then $T_{\mathrm{max}}>T_*>0$ thanks to \eqref{sep-loc1} and the continuity of $(u,v)$.
 Assume by contradiction, $T_{\mathrm{max}} < + \infty$. We observe that the \emph{a priori} estimates obtained in Section \ref{apro-es} hold for $(u,v)$ on $[0,T_{\mathrm{max}})$. In particular, the uniform strict separation property obtained in Proposition \ref{strictseparationvs} yields
\begin{align}
\| u(t)\|_{C(\overline{\Omega})}\leq 1- \delta_1,\quad
\| v(t) \|_{C(\overline{\Omega})}\leq 1- \delta_1, \quad  \forall \, t \in [0,T_{\mathrm{max}}).
\label{sep-loc2}
\end{align}
  As a consequence, $(u,v)$ is also a solution to the auxiliary problem \eqref{eq:as1} with the nonlinear term $F_{\delta_1/2}(u,v)$, that is
   $(u^\sharp,v^\sharp)=(u,v)$ on $[0,T_{\mathrm{max}})$.
   Since $(u^\sharp,v^\sharp)$ is global, by its continuity and \eqref{sep-loc2}, we have
   $$
   \| u^\sharp(T_{\mathrm{max}})\|_{C(\overline{\Omega})}\leq 1- \delta_1,\quad
   \| v^\sharp(T_{\mathrm{max}})\|_{C(\overline{\Omega})}\leq 1- \delta_1.
   $$
   Then there exists certain $t_0>0$ such that
   $$
   \| u^\sharp(t)\|_{C(\overline{\Omega})}\leq 1- \frac{3}{4}\delta_1,\quad
   \| v^\sharp(t)\|_{C(\overline{\Omega})}\leq 1- \frac{3}{4}\delta_1,\quad
   \forall\, t\in [T_{\mathrm{max}},T_{\mathrm{max}}+t_0],
   $$
   which imply $(u,v)=(u^\sharp,v^\sharp)$ is a strong solution to the original problem \eqref{eq:vs1}--\eqref{eq:vsi} on $[0,T_{\mathrm{max}}+t_0]$. This leads to a contradiction with the definition of $T_{\mathrm{max}}$.

   Hence, we can conclude $T_{\mathrm{max}}=+\infty$ and moreover, the \emph{a priori} estimates in Section \ref{apro-es} apply. Thanks to the strict separation property of the solution $(u,v)$ on $[0,+\infty)$, it enjoys the same regularity as $(u^\sharp,v^\sharp)$ (see Definition \ref{ssad}). In this way, we establish the existence of a (unique) global strong solution to problem \eqref{eq:vs1}--\eqref{eq:vsi} that fulfills the expected regularity properties and \eqref{sep:rstr}.
\qed
\medskip

We proceed to study the weak solutions of problem \eqref{eq:vs1}--\eqref{eq:vsi}.
\medskip

\noindent\textbf{Proof of Proposition \ref{wsve}.}
The proof is based on a suitable approximation for the initial data (cf. e.g., \cite{GP2022,MZ2004}).
Suppose that $\left(u_0 ,v_0\right) \in V$ satisfies $F(u_0,v_0) \in L^1(\Omega)$ and $\overline{u_0},\overline{v_0} \in [-1+m,1-m]$ for some $m \in (0,1)$. According to \cite{GP2022}, there exists a sequence $\{(u_{0,n},v_{0,n})\}_{n=1}^{\infty} \subset (H^3(\Omega)\cap W) \times (H^3(\Omega)\cap W)$ satisfying
\begin{align*}
    & \|u_{0,n} \|_{L^{\infty}(\Omega)}\le 1- \frac{1}{n},\quad \|v_{0,n} \|_{L^{\infty}(\Omega)} \le 1- \frac{1}{n}, \quad  \forall\, n \in \mathbb{Z}^+, \\
     & (u_{0,n},v_{0,n}) \to (u_0,v_0) \quad  \text{in} \ V \ \ \text{as} \ n \to +\infty.
\end{align*}
We can extract a subsequence of $\{(u_{0,n},v_{0,n})\}_{n=1}^{\infty}$ (not relabelled for simplicity) such that
\begin{equation*}
    (u_{0,n},v_{0,n}) \to (u_0,v_0) \quad \text{a.e. in} \ \Omega \ \ \mathrm{as} \ n \to +\infty
\end{equation*}
By the Lebesgue dominated convergence theorem, we find
\begin{equation*}
    \Psi(u_{0,n},v_{0,n}) \to \Psi(u_0,v_0) \quad \mathrm{as} \ n \to +\infty.
\end{equation*}
Then there exists a sufficiently large integer $N$ such that for all $n\ge N$, it holds
\begin{align*}
    & \|u_{0,n}\|_V \le 1+ \|u_0\|_V, \quad \|v_{0,n}\|_V \le 1+ \|v_0\|_V, \\
    & |\overline{u_{0,n}}|,\ |\overline{v_{0,n}}| \le 1-\frac{m}{2}, \quad \Psi(u_{0,n},v_{0,n}) \le 1 + \Psi(u_0,v_0).
\end{align*}

Consider the following approximate problem
\begin{equation}
\left\{\
\begin{aligned}
     & \partial_t u_n = \Delta \mu_n  &&\qquad\quad   \mathrm{in} \ \Omega \times (0,+\infty), \\
     & \mu_n= \alpha \partial_t u_n -\Delta u_n +  \partial_u F(u_n,v_n) &&\qquad \quad  \mathrm{in} \ \Omega \times (0,+\infty), \\
     & \partial_t v_n + \sigma (v_n -c) = \Delta \varphi_n &&\qquad\quad   \mathrm{in} \ \Omega \times (0,+\infty),\\
     & \varphi_n = \alpha \partial_t v_n  - \Delta v_n + \partial_v F(u_n,v_n)  && \qquad\quad  \mathrm{in} \ \Omega \times (0,+\infty),\\
     &\partial_\mathbf{n} u_n = \partial_\mathbf{n} v_n = \partial_\mathbf{n} \mu_n = \partial_\mathbf{n} \varphi_n = 0 && \qquad\quad  \mathrm{on} \ \partial \Omega \times (0,+\infty), \\
     & (u_n,v_n)|_{t=0} = (u_{0,n},v_{0,n}) &&\qquad\quad   \mathrm{in} \ \mathrm{\Omega}.
\end{aligned}
\right.
\label{eq:appro}
\end{equation}
It follows from Proposition \ref{ssve} that, for every $n\in \mathbb{Z}^+$, problem \eqref{eq:appro} admits a unique global strong solution $(u_n,v_n)$ on $[0,+\infty)$. Moreover, thanks to the above construction of initial data and Lemma \ref{dev},  $(u_n,v_n,\mu_n,\varphi_n)$ are uniformly bounded (w.r.t. $n$) in the following sense
\begin{align*}
    & u_n,v_n \in   L^{\infty}(0,+\infty;V)\cap L^2_{\mathrm{uloc}}(0,+\infty;W) \cap H^1_{\mathrm{uloc}}(0,+\infty;V^*),\\
    & \sqrt{\alpha}\partial_t u_n,\ \sqrt{\alpha}\partial_t v_n \in L^2_{\mathrm{uloc}}(0,+\infty;H),\\
    & \mu_n,\ \varphi_n \in L^2_{\mathrm{uloc}}(0,+\infty;V).
\end{align*}
Besides, $\widehat{S}'_{(u)}(u_n),\ \widehat{S}'_{(v)}(v_n)  \in L^2_{\mathrm{uloc}}(0,\infty;H)$ are uniformly bounded as well. This implies that
$u_n,v_n\in  L^\infty(\Omega\times(0,+\infty))$ and
$|u_n(x,t)|,\,|v_n(x,t)| <1 $ a.e. in $\Omega \times (0,+\infty)$.
Then there exists some functions
$(u,v,\mu,\varphi)$ with the same regularity properties and a convergent subsequence $\{(u_n,v_n,\mu_n,\varphi_n)\}$ such that as $n\to +\infty$,
$$(u_n,v_n,\mu_n,\varphi_n)\to (u,v,\mu,\varphi)$$
weakly (or weakly-$*$) in the corresponding spaces. Hereafter, the related convergence will always be understood in the sense of a subsequence.

By the Aubin-Lions-Simon lemma \cite{Si1987}, for any $T>0$, we obtain
$$u_n\to u, \ \ v_n \to v\quad\text{strongly in}\ \  C([0,T];H),$$
 which also implies
$$
u_n\to u,\ \ v_n \to v \quad\text{a.e. in}\ \ \Omega\times(0,T).
$$
As a consequence,
$$
\widehat{S}'_{(u)}(u_n)\to \widehat{S}'_{(u)}(u),\ \
\widehat{S}'_{(v)}(v_n)\to \widehat{S}'_{(v)}(v)\quad \text{weakly in}\ \ L^2(0,T;H),
$$
From $\mathbf{(H1)}$, we further infer that $u,\,v\in  L^\infty(\Omega\times(0,T))$ with $|u(x,t)|,\,|v(x,t)| <1 $ almost everywhere in $\Omega \times (0,T)$. Using $\mathbf{(H3)}$, we also find
\begin{equation*}
    \partial_u W(u_n,v_n)\to \partial_u W(u,v), \ \ \partial_v W(u_n,v_n)\to \partial_v W(u,v)\quad \text{weakly in}\ \ \ L^2(0,T;H).
\end{equation*}
Because of the strong convergence of $u_n\to u$, $v_n \to v$ in $C([0,T];H)$, we can verify the initial condition $(u,v)|_{t=0} = (u_0,v_0)$. Hence, it is straightforward to check that the limit $(u,v,\mu,\varphi)$ is a global weak solution to problem \eqref{eq:vs1}--\eqref{eq:vsi}. Uniqueness of the weak solution follows from Proposition \ref{ssvc}.

Next, we study the regularity of weak solutions. From Lemma \ref{hoev} and the construction of the initial data $\{(u_{0,n},v_{0,n})\}_{n=1}^{\infty}$, we find that for any $\kappa\in(0,1]$,
\begin{align*}
& \partial_t u_n,\ \partial_t v_n \in L^2_{\mathrm{uloc}}(\kappa,+\infty;V),\quad  \sqrt{\alpha}\partial_t u_n,\ \sqrt{\alpha}\partial_t v_n \in L^{\infty}(\kappa,+\infty;H),\\
& \mu_n,\varphi_n \in L^{\infty}(\kappa,+\infty;V),\quad  u_n,v_n \in L^{\infty}(\kappa,+\infty;W),\\
& \widehat{S}'_{(u)}(u_n),\ \widehat{S}'_{(v)}(v_n) \in L^{\infty}(\kappa,+\infty;H),
\end{align*}
with uniform bounds with respect to $n$ in the corresponding spaces. Then on the interval  $[\kappa,+\infty)$, the convergent subsequence $(u_n,v_n,\mu_n,\varphi_n)$ considered in (1) will converge (maybe up to a further subsequence) to some functions $(\widehat{u},\widehat{v},\widehat{\mu},\widehat{\varphi})$ with higher-order regularity properties described above. By uniqueness of the weak convergence, we have $(\widehat{u}(t),\widehat{v}(t),\widehat{\mu}(t),\widehat{\varphi}(t))=(u(t),v(t),\mu(t),\varphi(t))$ for all $t\geq \kappa$. Since $\kappa\in (0,1]$ is arbitrary, we see that every weak solution to problem \eqref{eq:vs1}--\eqref{eq:vsi} regularizes instantaneously as long as $t>0$. The Aubin-Lions-Simon lemma yields $u, v\in C([\kappa, \infty);H^{2-\epsilon}(\Omega))$ for any $\epsilon\in (0,1/2)$, so that $u, v\in C([\kappa, \infty);C(\overline{\Omega}))$ thanks to the Sobolev embedding theorem for $d\in\{2,3\}$. Then the strict separation property \eqref{rsep-w} is a consequence of Proposition \ref{strictseparationvs}.
\qed


\section{Regularity of Weak Solutions}\label{nv}
\setcounter{equation}{0}
In this section, we study the original problem \eqref{eq:nvs1}--\eqref{eq:nvsi}.
First, we recover the existence of global weak solutions by using an approach different from that in \cite{DG2022}.
\medskip

\noindent\textbf{Proof of Proposition \ref{hoenv}.}
Consider the regularized problem \eqref{eq:vs1}--\eqref{eq:vsi} subject to the given initial data $(u_0, v_0)$. For every $\alpha\in (0,1)$, thanks to Proposition \ref{wsve}, the problem admits a unique global weak solution, which we denote by
$(u_{\alpha},v_{\alpha},\mu_{\alpha},\varphi_{\alpha})$. Moreover, we have
\begin{align*}
& u_{\alpha}, v_{\alpha} \in C([0,+\infty);V)
\cap L^2_{\mathrm{uloc}}(0,+\infty;W) \cap H^1_{\mathrm{uloc}}(0,+\infty;V^*),\\
& \sqrt{\alpha}\partial_t u_{\alpha},\ \sqrt{\alpha}\partial_t v_{\alpha} \in L^2_{\mathrm{uloc}}(0,+\infty;H),\\
&\mu_{\alpha}, \varphi_{\alpha} \in  L^2_{\mathrm{uloc}}(0,+\infty;V),\quad \widehat{S}'_{(u)}(u_{\alpha}),\ \widehat{S}'_{(v)}(v_{\alpha})  \in L^2_{\mathrm{uloc}}(0,+\infty;H),
\end{align*}
with uniform bounds with respect to $\alpha\in (0,1)$ in the corresponding spaces.
In analogy to the proof of Proposition \ref{wsve}, we can find some
$(u,v,\mu,\varphi)$ with the same regularity properties and a convergent subsequence  $(u_{\alpha},v_{\alpha},\mu_{\alpha},\varphi_{\alpha})$ (nor relabelled for simplicity) such that
$$(u_{\alpha},v_{\alpha},\mu_{\alpha},\varphi_{\alpha})\to (u,v,\mu,\varphi) \quad \text{as}\ \alpha \to 0,$$
weakly (or weakly-$*$) in the corresponding spaces. In particular, for any $T>0$, it holds
\begin{equation*}
    \alpha \partial_t u_{\alpha}\to 0,\quad \alpha \partial_t v_{\alpha} \to 0\quad
     \text{strongly in}\ \ L^2(0,T;H).
\end{equation*}
By a standard compactness argument similar to that for Proposition \ref{wsve}, we can verify that $(u,v,\mu,\varphi)$ is a global weak solution to problem \eqref{eq:nvs1}--\eqref{eq:nvsi} on $[0,+\infty)$ with the initial data $(u_0,v_0)$. The uniqueness is a direct consequence of the continuous dependence estimate \eqref{es-conti1}, which can be proved by exactly the same argument as in \cite[Section 4.4]{DG2022}. Moreover, thanks to the regularity of weak solutions, the dissipative estimates \eqref{eq:denvs-1}, \eqref{eq:denvs-2} can be derived as in Lemma \ref{dev} with $\alpha=0$.
\qed
\medskip

Next, we show the regularizing effect of weak solutions for $t>0$.

\medskip

\noindent\textbf{Proof of Theorem \ref{hoenv}.}
Thanks to Lemma \ref{hoev}, for any $\kappa\in (0,1]$, we find that the convergent subsequence $\{(u_{\alpha},v_{\alpha},\mu_{\alpha},\varphi_{\alpha})\}$ in the proof of Proposition \ref{hoenv} satisfies
\begin{align*}
& u_{\alpha},v_{\alpha} \in L^{\infty}(\kappa,+\infty;W),\quad \partial_t u_{\alpha},\ \partial_t v_{\alpha} \in L^2_{\mathrm{uloc}}(\kappa,+\infty;V)\cap L^\infty(\kappa,+\infty;V^*),\\
& \mu_{\alpha},\varphi_{\alpha} \in L^{\infty}(\kappa,+\infty;V),\quad  \widehat{S}'_{(u)}(u_{\alpha}),\ \widehat{S}'_{(v)}(u_{\alpha}) \in L^{\infty}(\kappa,+\infty;H),
\end{align*}
with uniform bounds with respect to $\alpha\in (0,1)$ in the corresponding spaces. Based on the above facts, by the same argument as for Proposition \ref{wsve}, we can show that the weak solution $(u,v,\mu,\varphi)$ fulfills
\begin{align*}
& u,v \in L^{\infty}(\kappa,+\infty;W),\quad \partial_t u,\ \partial_t v \in L^2_{\mathrm{uloc}}(\kappa,+\infty;V)\cap L^\infty(\kappa,+\infty;V^*),\\
& \mu,\varphi \in L^{\infty}(\kappa,+\infty;V),\quad  \widehat{S}'_{(u)}(u),\ \widehat{S}'_{(v)}(u) \in L^{\infty}(\kappa,+\infty;H),
\end{align*}
with bounds in the corresponding spaces depending on ${\Psi}(u_0,v_0)$, $\overline{u_0}$, $\overline{v_0}$, $\Omega$, parameters of the system and $\kappa$.
By comparison in the equations \eqref{eq:nvs1}, \eqref{eq:nvs3} that hold almost everywhere in $\Omega\times [\kappa,+\infty)$, we get
$\mu, \varphi \in L^2_{\mathrm{uloc}}(\kappa,+\infty;H^3(\Omega))$.
Next, let us write the equations \eqref{eq:nvs2} and \eqref{eq:nvs4} into the following form:
\begin{align}
    & -\Delta u + \widehat{S}'_{(u)}(u) = \mu + \theta_{0,u}u - \partial_u W(u,v),\label{eq:ellipticu}\\
    & -\Delta v + \widehat{S}'_{(v)}(v) = \varphi + \theta_{0,v}v -\partial_v W(u,v). \label{eq:ellipticv}
\end{align}
The two terms in right-hand side of \eqref{eq:ellipticu} and \eqref{eq:ellipticv} are bounded in $L^{\infty}(\kappa,+\infty;V)$. Recalling that $\partial_\mathbf{n} u= \partial_\mathbf{n} v=0$ almost everywhere on $\partial\Omega\times [\kappa, +\infty)$, then we can apply Lemma \ref{nonlinearelliptic} to conclude
\begin{align*}
    &\| u \|_{L^\infty(\kappa,t;W^{2,p}(\Omega))}+ \| v \|_{L^\infty(\kappa,t;W^{2,p}(\Omega))} +
    \|\widehat{S}'_{(u)}(u)\|_{L^\infty(\kappa, t; L^p(\Omega))}\\
    &\quad + \|\widehat{S}'_{(v)}(v)\|_{L^\infty(\kappa, t; L^p(\Omega))}
    \leq C,\quad \forall\, t\geq \kappa,
\end{align*}
where $p = 6$ if $d=3$, or $p \in [2,+\infty)$ if $d=2$,
the positive constant $C$ depends on ${\Psi}(u_0,v_0)$, $\overline{u_0}$, $\overline{v_0}$, $\Omega$, parameters of the system and $\kappa$, but not on $t$.
\qed

\medskip

Finally, we prove the instantaneous strict separation property of $(u, v)$ in two dimensions.

\medskip

\noindent\textbf{Proof of Theorem \ref{strictseparation2d}}.
The proof essentially follows the idea in \cite{GP2023}. For any fixed $\kappa\in (0,1]$, let us consider the elliptic equation \eqref{eq:ellipticu} in $\Omega \times [\kappa,+\infty)$ subject to the boundary condition $\partial_\mathbf{n} u=0$ on $\partial\Omega \times [\kappa,+\infty)$.
As in the proof of \cite[Theorem 3.3]{GP2023}, we introduce the sequence:
\begin{equation*}
    k_n = 1-\delta_u -\frac{\delta_u}{2^n}, \quad \forall\, n\in \mathbb{Z}^+,
\end{equation*}
where $\delta_u\in (0,1)$ is a constant to be determined later.
Then it holds
\begin{equation*}
    1-2\delta_u < k_n < k_{n+1} < 1-\delta_u, \quad  \forall\, n \ge 1, \quad  k_n \to 1-\delta_u \ \ \mathrm{as} \ n \to +\infty.
\end{equation*}
For every $n\in \mathbb{Z}^+$, define
\begin{align*}
    & A_n(t) = \left\{x \in \Omega\ |\ u(x,t) - k_n \ge 0 \right\},
    \quad  \forall\, t\in [\kappa,+\infty),\\
    & u_n(x,t) = (u(x,t)-k_n)^+,\qquad z_n(t) = \int_{A_n(t)} 1\, \mathrm{d}x.
\end{align*}
Recalling that $\mu\in L^\infty(\kappa,+\infty,V)$, below we consider an arbitrary but fixed time $t_*\in [\kappa,+\infty)$ such that $\|\mu(t_*)\|_V\leq \|\mu\|_{L^\infty(\kappa,+\infty,V)}$. Testing the equation \eqref{eq:ellipticu} at $t_*$ with $u_n(x,t_*)$ yields
\begin{equation}
    \| \nabla u_n \|^2 + \int_{\Omega} \widehat{S}_{(u)}'(u) u_n\, \mathrm{d}x = \int_{\Omega} \left[\theta_{0,u}u -\partial_uW(u,v)\right] u_n\, \mathrm{d}x
    + \int_{\Omega} \mu u_n\, \mathrm{d}x.
    \label{sep-es1}
\end{equation}
By Theorem \ref{hoenv} and the Aubin-Lions-Simon lemma, we have $u, v\in C([\kappa, +\infty);H^{2-\epsilon}(\Omega))$ for any $\epsilon\in (0,1/2)$, so that $u, v\in C([\kappa, +\infty);C(\overline{\Omega}))$ thanks to the Sobolev embedding theorem for $d\in\{2,3\}$. Then it holds $|u(x,t_*)|\leq 1$ and $|v(x,t_*)|\le 1$ in $\Omega$, which imply
\begin{equation*}
    \int_{\Omega} \left[\theta_{0,u}u -\partial_uW(u,v)\right] u_n\, \mathrm{d}x \le \left( \max_{[-1,1]^2}|\partial_u W(\cdot,\cdot)|+ \theta_{0,u}\right)
    \int_{\Omega} u_n \, \mathrm{d}x,
\end{equation*}
and $0\leq u_n\leq 2\delta_u$. The other terms in \eqref{sep-es1} can be estimated exactly as in \cite[Section 3]{GP2023}, thus we get
\begin{align}
 & \| \nabla u_n \|^2 + \left(\widehat{S}_{(u)}'(1-2\delta_u)- \max_{[-1,1]^2}|\partial_u W(\cdot,\cdot)|-  \theta_{0,u}\right)\int_\Omega u_n\,\mathrm{d}x + \theta_u \int_\Omega u_n^2\,\mathrm{d}x \notag \\
  &\quad  =  \int_{\Omega} \mu u_n\, \mathrm{d}x = \int_{A_n(t_*)} \mu u_n\, \mathrm{d}x \notag \\
  &\quad  \leq 2\delta_u \|\mu(t_*)\|_{L^p(\Omega)}z_n^{1-\frac{1}{p}}
  \leq C\delta_u \|\mu(t_*)\|_{V}  \sqrt{p}z_n^{1-\frac{1}{p}},\quad \forall\ p\in [2,+\infty).
    \label{sep-es2}
\end{align}
In the last step, we have used the following Gagliardo-Nirenberg-type inequality in two dimensions (see \cite{B2010})
$$
\|f\|_{L^p(\Omega)}\leq C\sqrt{p}\|f\|^\frac{2}{p}\|f\|_V^{1-\frac{2}{p}},\quad \forall\, f\in V,\quad p\in [2,+\infty).
$$
With the inequality \eqref{sep-es2} and the additional assumption $\mathbf{(H2)}$, we can follow the same De Giorgi's iteration scheme as in \cite[Section 3]{GP2023} to show that there exists some sufficiently small $\delta_u\in (0,1)$, it holds $z_n(t_*)\to 0$ as $n\to +\infty$. Observe that $z_n(t_*)\to |\{x\in \Omega: u(x,t_*)\geq 1-\delta_u\}$, we find $\|(u(t_*)-(1-\delta_u))^+\|_{L^\infty(\Omega)}=0$. In a similar way, one can prove $\|(u(t_*)-(-1+\delta_u'))^-\|_{L^\infty(\Omega)}=0$ for some small $\delta'_u\in (0,1)$. Since $t_*\geq \kappa$ is arbitrary, from the space-time continuity of $u$, we get
$$\|u(t)\|_{C(\overline{\Omega})} \le 1 - \widetilde{\delta}_u,\quad \forall\,t\geq \kappa,\quad \text{with}\ \ \widetilde{\delta}_u=\min\left\{\delta_u, \delta_u'\right\}.$$
By the same argument, we can find some sufficiently small $\widetilde{\delta}_v\in (0,1)$ such that
$$\| v(t)\|_{C(\overline{\Omega})} \le 1 - \widetilde{\delta}_v,\quad \forall\,t\geq \kappa.$$
Taking $\delta_\kappa=\min\{\widetilde{\delta}_u, \widetilde{\delta}_v\}$, we arrive at the conclusion \eqref{isp-2d}.
\qed

\section{Long-time behavior}\label{longt}
\setcounter{equation}{0}
In this section, we study the long-time behavior of global weak solutions to problem \eqref{eq:nvs1}--\eqref{eq:nvsi}.

\subsection{Characterization of the $\omega$-limit set}
For any initial data $u_0, v_0 \in V$ with $F(u_0,v_0) \in L^1(\Omega)$ and $\overline{u_0},\overline{v_0} \in (-1,1)$, let $(u,v)$ be the corresponding global weak solution to problem \eqref{eq:nvs1}--\eqref{eq:nvsi} given by Proposition \ref{cdivnv}. First, it follows that
\begin{align}
\overline{u}(t)=\overline{u_0},\quad \overline{v}(t)-c= \left(\overline{v_0}-c\right)e^{-\sigma t},\quad \forall\,t\geq 0.
\label{mass-decay}
\end{align}
Next, we write problem \eqref{eq:nvs1}--\eqref{eq:nvsi} into the following equivalent form
\begin{align}
     & \partial_t u =\Delta \mu && \mathrm{in} \ \Omega \times (0,T),\label{eq:nvs1e}\\
     & \mu = -\Delta u + \partial_u F(u,v)&& \mathrm{in} \ \Omega \times (0,T), \label{eq:nvs2e}\\
     & \partial_t v + \sigma (\overline{v}-c)=\Delta \widetilde{\varphi} && \mathrm{in} \ \Omega \times (0,T), \label{eq:nvs3e}\\
     & \widetilde{\varphi} = - \Delta v + \partial_v F(u,v) + \sigma \mathcal{N}(v -\overline{v}) && \mathrm{in} \ \Omega \times (0,T), \label{eq:nvs4e}\\
     &\partial_\mathbf{n} u=\partial_\mathbf{n} v =\partial_\mathbf{n} \mu =\partial_\mathbf{n} \widetilde{\varphi}  = 0 && \mathrm{on} \ \partial \Omega \times (0,T), \label{eq:nvsbe} \\
     & (u,v)|_{t=0} = (u_0,v_0) && \mathrm{in} \ \Omega. \label{eq:nvsie}
\end{align}
Through a standard argument and the chain rule (see e.g., \cite[Section 3]{GGM2017} for the single Cahn-Hilliard-Oono equation), we can work with the new form \eqref{eq:nvs1e}--\eqref{eq:nvsbe} and show that the weak solution $(u,v)$ satisfies the following energy equality
\begin{equation}\label{energyeqdif}
    \frac{\mathrm{d}}{\mathrm{d}t} \widetilde{\Psi}(u(t),v(t))
    + \|\nabla \mu(t)\|^2 + \|\nabla \widetilde{\varphi}(t)\|^2 + \sigma (\overline{v}(t)-c)\int_\Omega \varphi(t)\,\mathrm{d}x = 0,
\end{equation}
almost everywhere in $(0,+\infty)$, where
the modified free energy is given by (cf. \eqref{eq:energyfunctionalr})
\begin{equation}
     \widetilde{\Psi}(u,v) := \Psi(u,v) + \sigma \|v - \overline{v}\|^2_{V^*_0}.
     \notag
\end{equation}
By a similar argument for Lemma \ref{rlow} (cf. \eqref{eq:avepsi}), we find  $|\overline{\varphi}|\leq C(1+\|\nabla \varphi\|)$. This implies $\sigma (\overline{v}(t)-c)\int_\Omega \varphi(t)\,\mathrm{d}x\in L^1(0,1)$. Furthermore, Lemma \ref{hoev} yields
$\varphi\in L^\infty(1,+\infty;V)$, then from \eqref{mass-decay} we get
$$
\left|\sigma (\overline{v}(t)-c)\int_\Omega \varphi(t)\,\mathrm{d}x\right|
\leq K_1e^{-\sigma t},\quad \text{for a.a.}\  t\geq 1,
$$
where $K_1>0$ depends on ${\Psi}(u_0,v_0)$, $\overline{u_0}$, $\overline{v_0}$, $\Omega$, and parameters of the system, but not on $t$.
As a consequence, the  energy equality \eqref{energyeqdif} implies that the map $t\to \widetilde{\Psi}(u(t),v(t))$ is absolutely continuous
on $[0,+\infty)$. This fact combined with the assumptions $\mathbf{(H1)}$, $\mathbf{(H3)}$ and the facts $u, v\in C([0,+\infty);H)$ yields $u, v\in C([0, +\infty);V)$.

 Set
\begin{align}
E(t)=  \widetilde{\Psi}(u(t),v(t))+\frac{K_1}{\sigma}e^{-\sigma t}, \quad \forall\, t\geq 1.
\notag
\end{align}
It follows from \eqref{energyeqdif} that
\begin{align}
 \frac{\mathrm{d}}{\mathrm{d}t}E(t) + \|\nabla \mu(t)\|^2 + \|\nabla \widetilde{\varphi}(t)\|^2 \leq 0, \quad \text{for a.a.}\  t\geq 1.
 \label{diff-E}
\end{align}
Since $\widetilde{\Psi}(u,v)$ is bounded from below (cf. \eqref{low-bd1}), we see that $E(t)$ converges as $t \to +\infty$, i.e., there exists some $E_\infty\in \mathbb{R}$ such that
\begin{equation}
    \lim_{t \to +\infty} E(t) = E_\infty.
    \label{decayE}
\end{equation}
In view of \eqref{mass-decay}, we have the convergence of energy as well
\begin{equation}
    \lim_{t \to +\infty} \widetilde{\Psi}(u(t),v(t)) = E_\infty. \label{lim-ener}
\end{equation}

Below we study the $\omega$-limit set.
\bp\label{omegalim}
Let the assumptions in Proposition \ref{cdivnv} be satisfied. For any initial data $u_0, v_0 \in V$ with $F(u_0,v_0) \in L^1(\Omega)$ and $\overline{u_0},\overline{v_0} \in (-1,1)$, we define the corresponding $\omega$-limit set by
\begin{align*}
\omega(u_0,v_0)& =\big\{(u_\infty,v_\infty)\,|\, u_\infty,v_\infty \in H^{2-\epsilon}(\Omega),\  \overline{u_\infty}=\overline{u_0}, \ \overline{v_\infty}=c, \ \text{there exist}\ \{t_n\}\nearrow +\infty \\
&\qquad \qquad \qquad\ \  \text{such that}\ u(t_n)\to u_\infty, \ \ v(t_n)\to v_\infty\  \text{in}\ H^{2-\epsilon}(\Omega)\big\},
\end{align*}
where $\epsilon\in (0,1/2)$. Then $\omega(u_0,v_0)$ is non-empty, bounded in $W\times W$ and compact in  $H^{2-\epsilon}(\Omega)\times H^{2-\epsilon}(\Omega)$. Moreover, every $(\uinf,\vinf)\in \omega(u_0,v_0)$ is a strong solution to the stationary problem
\eqref{sta0}.
\ep
\begin{proof}
Since the weak solution $(u, v)$ is uniformly bounded in $L^\infty(1,+\infty; W\times W)$, we find that the set $\omega(u_0,v_0)$ is non-empty and compact in  $H^{2-\epsilon}(\Omega)\times H^{2-\epsilon}(\Omega)$ for any $\epsilon\in (0,1/2)$. Besides, $\omega(u_0,v_0)$ is a bounded set in $W\times W$.
Next, we show that $\omega(u_0,v_0)$ consists of steady states satisfying \eqref{sta0} by using the argument in \cite{HP2001}.
Let $(u_\infty, v_\infty)\in \omega(u_0,v_0)$. From \eqref{mass-decay}, we have
$\overline{u_\infty}=\overline{u_0}$, $\overline{v_\infty}=c$.
By definition, there exists an unbounded sequence $\{t_n\}$ such that $t_{n+1}\geq t_n+1$, and
$$
\lim_{n\to+\infty}\|(u(t_n),v(t_n))-(u_\infty,v_\infty)\|_{H^{2-\epsilon}(\Omega)\times H^{2-\epsilon}(\Omega)}=0,
$$
for some $\epsilon\in (0,1/2)$. Without loss of generality, we also assume $\{(u(t_n),v(t_n))\}$ weakly converge to $(u_\infty,v_\infty)$ in $W\times W$.
Integrating \eqref{diff-E} gives
\begin{align}
\int_1^{\infty}\left(\|\nabla \mu(t)\|^2 + \|\nabla \widetilde{\varphi}(t)\|^2\right)\,\mathrm{d}t\leq C.
\label{es-uv-u}
\end{align}
By comparison in \eqref{eq:nvs1e}, \eqref{eq:nvs3e} and using \eqref{mass-decay}, \eqref{es-uv-u}, we further obtain
\begin{align}
\int_1^{\infty}\left(\|\partial_t u(t)\|_{V^*}^2 + \|\partial_tv(t)\|_{V^*}^2\right)\,\mathrm{d}t \leq C.
\label{es-uv-t}
\end{align}
Like in \cite{HP2001}, we infer from \eqref{es-uv-t} that
\begin{align}
\int_0^{1}\left(\|\partial_t u(t_n+\tau)\|_{V^*}^2 + \|\partial_tv(t_n+\tau)\|_{V^*}^2\right)\,\mathrm{d}\tau\to 0, \quad \text{as}\ n\to +\infty,
\notag
\end{align}
which implies
$$
\|u(t_n+\tau_1)-u(t_n+\tau_2)\|_{V^*}\to 0,\quad \|v(t_n+\tau_1)-v(t_n+\tau_2)\|_{V^*}\to 0, \quad \text{uniformly for}\ \tau_1,\,\tau_2\in [0, 1].
$$
By precompactness of the trajectory $(u(t),v(t))$ in $H^{2-\epsilon}(\Omega)\times H^{2-\epsilon}(\Omega)$, we have for every $\tau \in [0,1]$
$$
\|u(t_n+\tau)-u_\infty\|_{H^{2-\epsilon}(\Omega)}\to 0, \quad
\|v(t_n+\tau)-v_\infty\|_{H^{2-\epsilon}(\Omega)}\to 0, \quad \text{as}\ n\to +\infty,
$$
which also yields the almost everywhere convergence in $\Omega$ (up to a subsequence). Recalling the uniform boundedness of $\partial_vF(u, v)$ in $L^2(t_n, t_n+1;H)$, we get
$$
\partial_vF(u(t_n+\tau),v(t_n+\tau))\to \partial_vF(u_\infty,v_\infty),\quad \text{weakly in}\ L^2(0, 1;H)\ \text{as}\ n\to +\infty.
$$
Thus, possibly up to a further subsequence, we deduce from \eqref{eq:nvs4e} that
$$
\widetilde{\varphi}(t_n+\tau)=\widetilde{\varphi}(u(t_n+\tau),v(t_n+\tau)) \to \widetilde{\varphi}_\infty,\quad \text{weakly in}\ L^2(0, 1;H)\ \text{as}\ n\to +\infty.
$$
where the limit can be identified as
\begin{align}
\widetilde{\varphi}_\infty=-\Delta v_\infty + \partial_vF(u_\infty,v_\infty) +\sigma \mathcal{N}(v_\infty-c)\in H.
\label{phi-inf}
\end{align}
We note that $\widetilde{\varphi}_\infty$ is independent of time.
Thus, for any $\eta\in H_0$, it holds
\begin{align}
(\widetilde{\varphi}_\infty, \eta)
& = \int_0^1 (\widetilde{\varphi}_\infty, \eta) \,\mathrm{d}\tau  = \lim_{n\to +\infty} \int_0^1 (\widetilde{\varphi}(t_n+\tau), \eta) \,\mathrm{d}\tau. \notag
\end{align}
It follows from \eqref{es-uv-u} and the Poincar\'{e}-Wirtinger inequality that
\begin{align}
& \left|\int_0^1 (\widetilde{\varphi}(t_n+\tau), \eta) \,\mathrm{d}\tau\right|  \leq C \int_0^1 \|\nabla \widetilde{\varphi}(t_n+\tau)\|\|\eta\|\,\mathrm{d}\tau\to 0, \quad \text{as}\ n\to +\infty.
\notag
\end{align}
Therefore, $(\widetilde{\varphi}_\infty, \eta)=0$ for any $\eta\in H_0$,
which implies that $\widetilde{\varphi}_\infty$ is a constant. Integrating \eqref{phi-inf} over $\Omega$ easily gives $\widetilde{\varphi}_\infty=\overline{\partial_vF(u_\infty,v_\infty)}$.
By a similar argument, we obtain
$$\mu_\infty=-\Delta \uinf + \partial_u F(\uinf,\vinf)\quad \text{with}\ \ \mu_\infty = \overline{\partial_u F(\uinf,\vinf)}.
$$
Hence, $(u_\infty,v_\infty)\in W\times W$ is a strong solution to the stationary problem
\eqref{sta0}.
\end{proof}

Given an initial datum $(u_0,v_0)$ as that in Proposition \ref{omegalim}, we denote its corresponding set of steady states by
$$
\mathcal{S}(u_0,v_0)=\left\{\,(u_{\mathrm{s}}, v_{\mathrm{s}})\in W\times W\,|\, (u_{\mathrm{s}}, v_{\mathrm{s}}) \ \ \text{a strong solution to}\ \ \eqref{sta0}\,\right\}.
$$
Proposition \ref{omegalim} implies that $\omega(u_0,v_0)\subset \mathcal{S}(u_0,v_0)$ and thus $\mathcal{S}(u_0,v_0)$ is non-empty. From \eqref{sta0}, we see that $\mathcal{S}(u_0,v_0)$ only relates the initial data via the mean value $\overline{u_0}$.
Moreover, we have

\bp[Strict separation property of steady states] \label{separationequilibrium}
Let the assumptions of Proposition \ref{omegalim} be satisfied.

(1) For every   $(u_{\mathrm{s}}, v_{\mathrm{s}}) \in \mathcal{S}(u_0,v_0)$, there exists some constant $\delta_{\mathrm{s}} \in (0,1)$ such that
\begin{equation*}
    \|u_{\mathrm{s}}\|_{C(\overline{\Omega})}\leq 1-\delta_{\mathrm{s}}, \quad \|v_{\mathrm{s}}\|_{C(\overline{\Omega})}\leq 1-\delta_{\mathrm{s}}.
\end{equation*}

(2) For all $(u_{\infty},v_{\infty})\in \omega(u_0,v_0)$, it holds
\begin{equation*}
    \|u_{\infty}\|_{C(\overline{\Omega})}\leq 1-\delta_{\infty}, \quad \|v_{\infty}\|_{C(\overline{\Omega})}\leq 1-\delta_{\infty},
\end{equation*}
where the constant $\delta_{\infty} \in (0,1)$ is independent of $(u_{\infty},v_{\infty})$.
\ep

\begin{proof}
We apply a dynamic approach inspired by \cite{MZ2004} (see its \cite{H2022} for further application), for possible alternative proof involving the maximum principle, we refer to \cite[Proposition 6.1]{AW2007} for the single Cahn-Hilliard equation. It is obvious that every steady state $(u_{\mathrm{s}}, v_{\mathrm{s}})\in \mathcal{S}(u_0,v_0)$ can be viewed as (at least) a global weak solution to the viscous problem
\eqref{eq:vs1}--\eqref{eq:vsb} with some $\alpha\in (0,1)$ and the initial data given by $(u_{\mathrm{s}}, v_{\mathrm{s}})$ itself. Then by the uniqueness of solutions and Proposition \ref{wsve}, we find there exists some $\delta_{\mathrm{s}}\in(0,1)$ such that
\begin{equation*}
    \|u_{\mathrm{s}}(t)\|_{C(\overline{\Omega})}\leq 1-\delta_{\mathrm{s}}, \quad \|v_{\mathrm{s}}(t)\|_{C(\overline{\Omega})}\leq 1-\delta_{\mathrm{s}}, \quad  \forall\, t \ge 1.
\end{equation*}
Since $(u_{\mathrm{s}}, v_{\mathrm{s}})$ is actually independent of time, this gives our first conclusion. The second conclusion easily follows from the compactness of  $\omega(u_0,v_0)$ in $H^{2-\epsilon}(\Omega)\times H^{2-\epsilon}(\Omega)$, which is continuously embedded in $C(\overline{\Omega}) \times C(\overline{\Omega})$ for $\epsilon\in (0,1/2)$, cf. \cite{AW2007}.
\end{proof}

Now we are able to show the eventual strict separation property of weak solutions.

\noindent\textbf{Proof of Theorem \ref{longtimesep}.}
It follows from relative compactness of the trajectory $(u(t), v(t))$ and Proposition \ref{omegalim} that
\begin{equation}\label{h2-sconvergence}
    \lim_{t \to +\infty} \mathrm{dist}((u(t),v(t)),\omega(u_0,v_0)) = 0 \quad  \mathrm{in} \  H^{2-\epsilon}(\Omega)\times H^{2-\epsilon}(\Omega).
\end{equation}
This combined with Proposition \ref{separationequilibrium}-(2) yields the conclusion \eqref{esp-3d}
\qed

\subsection{An extended {\L}ojasiewicz-Simon inequality}
The eventual strict separation property Theorem \ref{longtimesep} plays an essential role in the study of long-time behavior of global weak solutions (see \cite{DG2022} for the conserved case). Since we are only interested in the behavior of solutions as $t\to +\infty$, Theorem \ref{longtimesep} implies that the singularities of $\widehat{S}_{(u)}$, $\widehat{S}_{(v)}$ and their derivatives can be avoided by altering these functions outside a certain compact interval  $[-1+\delta,1-\delta]$ for some $\delta\in (0,\delta_{\mathrm{SP}})$, see \cite{AW2007,A2009, M2019} and the references therein.


Using a similar idea, under the additional assumption $\mathbf{(H4)}$ on the analyticity of $\widehat{S}_{(u)}$, $\widehat{S}_{(v)}$ and $W(u,v)$, the authors of \cite{DG2022} proved a suitable {\L}ojasiewicz-Simon type inequality in the conserved case. Below we present an extended version that works for the off-critical case, i.e., the change of mass is allowed.

\bp[{\L}ojasiewicz-Simon inequality] \label{ls}
Suppose that $\Omega \subset \mathbb{R}^d$ $(d\in\{2,3\})$ is a bounded domain with smooth boundary $\partial \Omega$, and the assumptions $(\mathbf{H0})$, $(\mathbf{H1})$, $(\mathbf{H3})$, $(\mathbf{H4})$ are satisfied. Given an initial data
$u_0, v_0 \in V$ with $F(u_0,v_0) \in L^1(\Omega)$ and $\overline{u_0},\overline{v_0} \in (-1,1)$, then for any $(u_\infty,v_\infty) \in \omega(u_0,v_0)$, there exist constants $\theta \in (0,1/2)$ and $C,\varpi > 0$ such that
for $u,v\in V\cap L^\infty(\Omega)$ satisfying $\norm{u}_{\lp{\infty}}\leq 1 -3\delta_\infty/4$, $\norm{v}_{\lp{\infty}} \leq 1 -3\delta_\infty/4$,   $\norm{u-\uinf}_{V} \leq \varpi$, $\norm{v-\vinf}_{V} \leq \varpi$, the following inequality holds
\begin{equation}\label{lsineq}
    \big|\widetilde{\Psi}(u,v)-\widetilde{\Psi}(u_\infty,v_\infty)\big|^{1-\theta}
    \le C \left( \|\mu-\overline{\mu}\|_{V^*_0}+ \|\widetilde{\varphi}-\overline{\widetilde{\varphi}}\|_{V^*_0} +  \left(|\overline{u}-\overline{u_0}| + |\overline{v}-c|\right)^{1-\theta}\right),
\end{equation}
where $\mu=-\Delta u+ \partial_uF(u,v)$, $\widetilde{\varphi} = - \Delta v + \partial_v F(u,v) + \sigma \mathcal{N}(v -\overline{v})$ and $\delta_\infty\in (0,1)$ is the constant determined in Proposition \ref{separationequilibrium}-(2).
\ep

\begin{proof}
\noindent\textbf{Step 1.}
In view of Proposition  \ref{separationequilibrium}-(2), we define
$Q=[-1+\delta_\infty/4,1-\delta_\infty/4]^2$ and
$$
F_{\mathrm{reg}}(s_1,s_2)= F(s_1,s_2)\chi_Q(s_1,s_2)+G(s_1,s_2)(1-\chi_Q(s_1,s_2)),
$$
where $G(s_1,s_2)$ is chosen in such a way to extend $F$ outside $Q$ with $C^3(\mathbb{R}^2)$ regularity and bounded derivatives up to order three (cf.  \cite{DG2022}).
Consider the regularized energy functional
$$
\widetilde{\Psi}_{\mathrm{reg}}(u,v) := \int_{\Omega}\left[ \frac{1}{2} |\nabla u|^2
    + \frac{1}{2} |\nabla v|^2 + F_{\mathrm{reg}}(u,v) \right] \mathrm{d}x + \sigma \|v - \overline{v}\|^2_{V^*_0}.
$$
It has been shown that $\widetilde{\Psi}_{\mathrm{reg}}(u,v)$ is twice continuously Fr\'{e}chet differentiable and $(u_\infty, v_\infty)$ is its critical point (see \cite[Lemma 7.3, Lemma 7.4]{DG2022}). Recall the first Fr\'{e}chet derivative 
$$
\left\langle \widetilde{\Psi}'_{\mathrm{reg}}(u,v), (h,k)\right\rangle
=\int_\Omega \big(\nabla u\cdot\nabla h + \nabla v\cdot\nabla k+ \partial_uF(u,v)h +\partial_vF(u,v)k +\sigma\mathcal{N}(v-\overline{v})k\big)\,\mathrm{d}x,
$$
for any $h, k\in V_0$. Under the additional assumption $\mathbf{(H4)}$, the following {\L}ojasiewicz--Simon type inequality has been obtained in \cite[Proposition 7.2]{DG2022}: there exist
$\theta_1 \in (0,1/2)$ and $C_1,\varpi_1> 0$ such that
for $\widehat{u},\widehat{v}\in V$ satisfying $\overline{\widehat{u}}= \overline{u_0}$, $\overline{\widehat{v}}= c$ and
$\norm{\widehat{u}-\uinf}_{V} \leq \varpi_1$, $\norm{\widehat{v}-\vinf}_{V} \leq \varpi_1$, it holds
\begin{align}
    \big|\widetilde{\Psi}_{\mathrm{reg}}(\widehat{u},\widehat{v}) -\widetilde{\Psi}_{\mathrm{reg}}(u_\infty,v_\infty)\big|^{1-\theta_1}
    &\le C_1 \big\|\widetilde{\Psi}_{\mathrm{reg}}'(\widehat{u},\widehat{v})\big\|_{V^*_0\times V^*_0}.
    \label{lsineq-reg}
\end{align}
From the construction of $\widetilde{\Psi}_{\mathrm{reg}}(\widehat{u},\widehat{v})$, if we further require that
$\widehat{u},\widehat{v}\in L^\infty(\Omega)$ satisfying $\norm{\widehat{u}}_{\lp{\infty}}\leq 1 -\delta_\infty/2$, $\norm{\widehat{v}}_{\lp{\infty}} \leq 1 -\delta_\infty/2$, then in \eqref{lsineq-reg}, we can simply replace $\widetilde{\Psi}_{\mathrm{reg}}(\widehat{u},\widehat{v})$ by $\widetilde{\Psi}(\widehat{u},\widehat{v})$, that is
\begin{equation}\label{lsineq-sing}
    \big|\widetilde{\Psi}(\widehat{u},\widehat{v}) -\widetilde{\Psi}(u_\infty,v_\infty)\big|^{1-\theta_1}
    \le C_1 \big\|\widetilde{\Psi}'(\widehat{u},\widehat{v})\big\|_{V^*_0\times V^*_0}  \le C_1\left( \|\widehat{\mu}-\overline{\widehat{\mu}}\|_{V^*_0}+ \|\widehat{\varphi}-\overline{\widehat{\varphi}}\|_{V^*_0}\right),
\end{equation}
where
$$
\widehat{\mu}=-\Delta \widehat{u}+ \partial_uF(\widehat{u},\widehat{v}),\quad \widehat{\varphi} = - \Delta \widehat{v} + \partial_v F(\widehat{u},\widehat{v}) + \sigma \mathcal{N}(\widehat{v} -\overline{\widehat{v}}).
$$

 \textbf{Step 2.} We shall control the possible mass change by using the perturbation argument in \cite{HW2023} (see also \cite{CRW2019}). For any functions $u,v\in V\cap L^\infty(\Omega)$ satisfying $\norm{u}_{\lp{\infty}}\leq 1 -3\delta_\infty/4$, $\norm{v}_{\lp{\infty}} \leq 1 -3\delta_\infty/4$ and $\norm{u-\uinf}_{V} \leq \varpi_2$, $\norm{v-\vinf}_{V} \leq \varpi_2$, let us define 
 \begin{align*}
 &\bracket{\widehat{u},\widehat{v}} = (u,v) - \bracket{\overline{u}-\overline{u_0},\overline{v}-c},
 \\
 &\mu=-\Delta u+ \partial_uF(u,v),\quad \widetilde{\varphi} = - \Delta v + \partial_v F(u,v) + \sigma \mathcal{N}(v -\overline{v}). 
 \end{align*} 
 We can choose the constant $\varpi_2\in (0,\varpi_1/2)$ sufficiently small such that
\begin{equation*}
    |\overline{u}-\overline{u_0}|\leq \min\left\{\frac{\delta_\infty}{4},\ \frac{\varpi_1}{2\sqrt{\Omega}}\right\},\quad |\overline{v}-c|\leq
     \min\left\{ \frac{\delta_\infty}{4},\ \frac{\varpi_1}{2\sqrt{\Omega}}\right\}.
\end{equation*}
Then it holds $\norm{\widehat{u}}_{\lp{\infty}}\leq 1 -\delta_\infty/2$, $\norm{\widehat{v}}_{\lp{\infty}}  \leq 1 -\delta_\infty/2$,
$\norm{\widehat{u}-\uinf}_{V} \leq \varpi_1$, $\norm{\widehat{v}-\vinf}_{V} \leq \varpi_1$. As a consequence, $\bracket{\widehat{u},\widehat{v}}$ satisfies the inequality \eqref{lsineq-sing}.
On the other hand, since both $(u,v)$ and $\bracket{\widehat{u},\widehat{v}}$ are strictly separated from $\pm 1$, we infer from $\mathbf{(H1)}$ and $\mathbf{(H3)}$ that
\begin{align}
        \|(\mu-\overline{\mu})- (\widehat{\mu}-\overline{\widehat{\mu}})\|_{V_0^*} + 
        \|(\varphi-\overline{\varphi})- (\widehat{\varphi}-\overline{\widehat{\varphi}})\|_{V_0^*}
        &\le C \bracket{\abs{\overline{u}-\overline{u_0}} + \abs{\overline{v}-c}} \notag \\
        &\le C \bracket{\abs{\overline{u}-\overline{u_0}} + \abs{\overline{v}-c}}^{1-\theta_1},
    \label{e'difference}
\end{align}
as well as
\begin{equation}\label{edifference}
    \big|\widetilde{\Psi}(u,v) -\widetilde{\Psi}(\widehat{u},\widehat{v}) \big| ^{1-\theta_1}
    \le C \bracket{\abs{\overline{u}-\overline{u_0}} + \abs{\overline{v}-c}}^{1-\theta_1}.
\end{equation}
Inserting \eqref{edifference} and \eqref{e'difference} into \eqref{lsineq-sing}, we find
\begin{equation}
    \begin{aligned}
         \big|\widetilde{\Psi}(u,v) -\widetilde{\Psi}(u_\infty,v_\infty)\big|^{1-\theta_1}
        & \le \big|\widetilde{\Psi}(\widehat{u},\widehat{v}) -\widetilde{\Psi}(u_\infty,v_\infty)\big|^{1-\theta_1} + \big|\widetilde{\Psi}(u,v) -\widetilde{\Psi}(\widehat{u},\widehat{v}) \big| ^{1-\theta_1}  \\
        & \le C_1\left( \|\widehat{\mu}-\overline{\widehat{\mu}}\|_{V^*_0}+ \|\widehat{\varphi}-\overline{\widehat{\varphi}}\|_{V^*_0}\right) + C \bracket{\abs{\overline{u}-\overline{u_0}} + \abs{\overline{v}-c}}^{1-\theta_1}\\
        & \le C_1 \left(\|\mu-\overline{\mu}\|_{V^*_0}+ \|\widetilde{\varphi}-\overline{\widetilde{\varphi}}\|_{V^*_0} \right) + C_2 \bracket{\abs{\overline{u}-\overline{u_0}} + \abs{\overline{v}-c}}^{1-\theta_1}.\notag
    \end{aligned}
\end{equation}
Hence, taking $\theta=\theta_1$, $\varpi=\varpi_2$ and $C=\max\{C_1,C_2\}$, we arrive at the conclusion \eqref{lsineq}.
\end{proof}

\subsection{Convergence to equilibrium}

Now we are able to prove the convergence to a single equilibrium by using the well-known {\L}ojasiewicz-Simon approach. 
\smallskip 

\noindent\textbf{Proof of Theorem \ref{convergequi}.} With the aid of Proposition \ref{ls} and the auxiliary energy inequality \eqref{diff-E}, we can apply an argument similar to that in \cite{DG2022,GGM2017} with minor modification. The the convenience of the readers, we sketch it here.  

It follows from \eqref{lim-ener} that 
\begin{align}
 \widetilde{\Psi}(u_\infty,v_\infty)=E_\infty,\quad \forall\, (u_\infty,v_\infty)\in \omega(u_0,v_0).
 \notag
\end{align} 
Next, thanks to the compactness of $\omega(u_0,v_0)$ in $\hs{2-\epsilon}\times \hs{2-\epsilon}\subset V\times V$ $(\epsilon\in (0,1/2))$, we can find a finite number of steady states $\big(u_\infty^{(j)},v_\infty^{(j)}\big)\in \omega(u_0,v_0)$ $(1\leq j\leq N)$ such that
$$
\omega(u_0,v_0)\subset \mathcal{B}=\bigcup_{j=1}^N B\left((u_\infty^{(j)},v_\infty^{(j)}), \varpi_j\right),
$$ 
where $\varpi_j>0$ is the constant corresponding to $(u_\infty^{(j)},v_\infty^{(j)})$ in Proposition \ref{ls} and $B\big((u_\infty^{(j)},v_\infty^{(j)}), \varpi_j\big)$ denotes the open ball in $V\times V$, centered at $(u_\infty^{(j)},v_\infty^{(j)})$ with radius $\varpi_j$. By Proposition \ref{separationequilibrium}-(2) and \eqref{h2-sconvergence}, there exists a sufficiently large $T_1\geq 1 $ such that
\begin{equation*}
    \norm{u(t)}_{C(\overline{\Omega})} \leq 1 -\frac{3}{4}\delta_{\infty},\quad \norm{v(t)}_{C(\overline{\Omega})} \leq 1 -\frac{3}{4}\delta_{\infty},\quad 
    (u(t),v(t))\in \mathcal{B},\quad 
      \forall\, t\geq T_1.
\end{equation*}
By Proposition \ref{ls}, we can extract uniform constants $C > 0$, $\theta\in(0, 1/2)$ such that for all $t\geq T_1$ the following inequality holds 
\begin{equation}\label{lsineq-sol}
    \big|\widetilde{\Psi}(u(t),v(t))-E_\infty\big|^{1-\theta}
    \le C \left( \|\mu(t)-\overline{\mu}(t)\|_{V^*_0}+ \|\widetilde{\varphi}(t)-\overline{\widetilde{\varphi}}(t)\|_{V^*_0} +   |\overline{v}(t)-c|^{1-\theta}\right),
\end{equation}
where $\mu=-\Delta u+ \partial_uF(u,v)$, $\widetilde{\varphi} = - \Delta v + \partial_v F(u,v) + \sigma \mathcal{N}(v -\overline{v})$. 
Therefore, we infer from \eqref{mass-decay}, \eqref{diff-E}, \eqref{decayE} and \eqref{lsineq-sol} that for any $t\geq T_1$,
\begin{align}
        &\int_{t}^\infty \left(\norm{\nabla \mu(\tau)}^2 +\norm{\nabla \tvp(\tau)}^2\right)\,\mathrm{d}\tau \notag \\
        & \quad \leq  E(t)-E_\infty \leq |\widetilde{\Psi}(u(t),v(t)) -E_\infty|+\frac{K_1}{\sigma}e^{-\sigma t} \notag \\
        &\quad  \le C \left( \|\mu(t)-\overline{\mu}(t)\|_{V^*_0}^\frac{1}{1-\theta}+ \|\widetilde{\varphi}(t)-\overline{\widetilde{\varphi}}(t)\|_{V^*_0} ^\frac{1}{1-\theta}+ |\overline{v}(t)-c|\right) + \frac{K_1}{\sigma}e^{-\sigma t} \notag \\
        &\quad  \le C \left(\norm{\nabla \mu(t)}^{\frac{1}{1-\theta}} + \norm{\nabla \tvp(t)}^{\frac{1}{1-\theta}} +  e^{-\sigma t}\right).
        \label{con-ener}
\end{align}
Since $\theta \in \bracket{0,1/2}$, we have
\begin{equation*}
    \int_t^\infty e^{-2 \sigma (1-\theta)\tau}\,\mathrm{d}\tau = \frac{1}{2 \sigma (1-\theta)}e^{-2 \sigma (1-\theta)t}\leq \frac{1}{2 \sigma (1-\theta)}e^{- \sigma t},\quad  \forall t \ge 1.
\end{equation*}
Set $\Sigma(t) = \norm{\nabla \mu(t)}+\norm{\nabla \tvp(t)}+ e^{-\sigma(1-\theta)t}$, we infer from \eqref{con-ener} that 
\begin{equation}
    \int_t^\infty \Sigma(\tau)^2\,\mathrm{d}\tau \le C \Sigma(t)^{\frac{1}{1-\theta}}, \ \quad \forall t\ge T_1. \notag 
\end{equation}
It follows from the above estimate and \cite[Lemma 7.1]{FS2000} that $\Sigma \in L^1(T_2,+\infty)$. By comparison in the evolution
equations \eqref{eq:nvs1e}, \eqref{eq:nvs3e} and using \eqref{mass-decay} we also have $\partial_t u, \partial_t v \in L^1(T_2,\infty;V^*)$, which implies the convergence of $(u(t), v(t))$ in $V^*\times V^*$ to some $(\uinf,\vinf)\in \omega(u_0,v_0)$ as $t\to +\infty$. By compactness of the trajectory, we can conclude
\begin{equation}
    \lim_{t \to +\infty} \norm{(u(t),v(t))-(\uinf,\vinf)}_{\hs{2-\epsilon} \times \hs{2-\epsilon}} = 0 \quad  \mathrm{for \ all} \ \epsilon \in \bracket{0,\half}.
    \notag 
\end{equation}
As a consequence, $\omega(u_0,v_0)$ is a singleton. 
\qed
%

\section{Appendix} 
\setcounter{equation}{0}

\subsection{Well-posedness of the auxiliary problem with regularized potential}
In what follows, we analyze the auxiliary problem \eqref{eq:as1}. First, we give the definition of strong solutions.
\bd \label{ssad}
Let $d\in \{2,3\}$, $T \in (0,+\infty)$ and $\alpha, \delta\in (0,1)$. For any initial data
$u_0,v_0 \in W$, the quadruple $(u,v,\mu,\varphi)$ is called a strong solution to problem \eqref{eq:as1} on $[0,T]$, if
\begin{align*}
    & u,v \in C([0,T];W) \cap L^2(0,T;H^3(\Omega)), \\
    & \partial_t u, \ \partial_t v \in C([0,T];H) \cap L^2(0,T;V) \cap H^1(0,T;V^*),\\
    & \mu, \varphi \in C([0,T];W) \cap L^2(0,T;H^3(\Omega))\cap H^1(0,T;V),
\end{align*}
the equations in \eqref{eq:as1} are
satisfied almost everywhere in $\Omega \times (0,T)$ and the initial conditions are satisfied everywhere in $\Omega$.
\ed

Then we prove the following result:

\bp \label{ssae}
Suppose that $\Omega \subset \mathbb{R}^d$ $(d\in\{2,3\})$ is a bounded domain with smooth boundary $\partial \Omega$, $T\in (0,+\infty)$, $\alpha, \delta \in (0,1)$ and the assumptions $(\mathbf{H0})$, $(\mathbf{H1})$, $(\mathbf{H3})$ are satisfied.
Then for any initial data $u_0,v_0 \in W$, problem \eqref{eq:as1} admits a unique strong solution $(u,v,\mu,\varphi)$ on $[0,T]$ in the sense of Definition \ref{ssad}.
\ep
\begin{proof}
Let $g,h\in L^2(0,T;V)\cap H^1(0,T;V^*)$ be two given functions. We consider the following linear problem:
\begin{equation}\label{iterationsystem}
\left\{\
    \begin{aligned}
        & \partial_t u =\Delta \mu && \quad \mathrm{in} \ \Omega \times (0,T),\\
     & \mu = \alpha \partial_t u - \Delta u + g && \quad \mathrm{in} \ \Omega \times (0,T), \\
     & \partial_t v + \sigma (v-c)=\Delta \varphi && \quad \mathrm{in} \ \Omega \times (0,T), \\
     & \varphi = \alpha\partial_t v - \Delta v + h && \quad \mathrm{in} \ \Omega \times (0,T),\\
     &\partial_{\mathbf{n}}u  = \partial_{\mathbf{n}} v  = \partial_{\mathbf{n}} \mu  = \partial_{\mathbf{n}} \varphi = 0 && \quad \mathrm{on} \ \partial \Omega \times (0,T), \\
     & (u,v)|_{t=0} = (u_0,v_0) && \quad \mathrm{in} \ \Omega,
    \end{aligned}
    \right.
\end{equation}
where the initial data satisfy $u_0, v_0\in W$.
Set
$$
\widetilde{u}=u-\overline{u},\quad \widetilde{v}=v-\overline{v},\quad  \widetilde{g}=g-\overline{g},\quad \widetilde{h}=h-\overline{h},
$$
with
\begin{equation}\notag
\overline{u}(t) = \overline{u_0}\quad \text{and}\quad
\overline{v}(t) = \overline{v_0} e^{-\sigma t}+ c\left(1-e^{-\sigma t}\right),
\quad \forall\, t\in [0,T].
\end{equation}
Then problem \eqref{iterationsystem} can be written as (cf. \cite{MZ2005})
\begin{equation}\label{iterationsystem-a}
\left\{\
    \begin{aligned}
     & \alpha \partial_t \widetilde{u} = \Delta \widetilde{u} +(\alpha I+ \mathcal{N})^{-1}\widetilde{u} -\alpha (\alpha I+ \mathcal{N})^{-1} \widetilde{g} && \quad \mathrm{in} \ \Omega \times (0,T), \\
     & \alpha \partial_t \widetilde{v} = \Delta \widetilde{v} + (\alpha I+ \mathcal{N})^{-1}\widetilde{v}  - \alpha \sigma (\alpha I+ \mathcal{N})^{-1}\widetilde{v} -\alpha (\alpha I+ \mathcal{N})^{-1}\widetilde{h}  && \quad \mathrm{in} \ \Omega \times (0,T),\\
     &\partial_{\mathbf{n}}\widetilde{u}  = \partial_{\mathbf{n}} \widetilde{v}  = 0 && \quad \mathrm{on} \ \partial \Omega \times (0,T), \\
     & (\widetilde{u},\widetilde{v})|_{t=0} = (u_0-\overline{u_0},v_0-\overline{v_0}) && \quad \mathrm{in} \ \Omega,
    \end{aligned}
    \right.
\end{equation}
in which we find two linear heat equations with compact perturbations. By a standard argument (e.g., the Faedo-Garlerkin approximation), we can show that problem \eqref{iterationsystem-a} admits a unique strong solution $(\widetilde{u},\widetilde{v})$ satisfying
$$\widetilde{u}, \widetilde{v}\in C([0,T];W)\cap L^2(0,T;H^3(\Omega))\cap H^1(0,T;V)\cap H^2(0,T;V^*).$$
Hence, problem \eqref{iterationsystem} admits a unique strong solution $(u,v,\mu, \varphi)$ with corresponding regularity properties.

Now for any given functions $u_*, v_*\in L^2(0,T; W) \cap H^1(0,T;H)$ satisfying $u_*(0)=u_0$, $v_*(0)=v_0$, we set
$$g_*= \partial_u F_\delta(u_*,v_*),\quad h_*= \partial_v F_\delta(u_*,v_*).$$
Since $F_\delta(\cdot,\cdot)\in BUC^2(\mathbb{R}\times\mathbb{R})$, it is straightforward to check that $g_*, h_*\in L^2(0,T;V)\cap H^1(0,T;V^*)$, which also implies $g_*, h_*\in C(0,T;H)$ thanks to the Aubin-Lions-Simon lemma \cite{Si1987}.
As a consequence, the mapping
\begin{equation*}
  \mathcal{M}:\ L^2(0,T; W)^2 \cap H^1(0,T;H)^2 \to L^2(0,T;H^3)^2\cap H^1(0,T;V)^2\ \ \text{s.t.}\ \   \mathcal{M}(u_*, v_*)=(u,v),
\end{equation*}
is well-defined, where $(u,v)$ is the unique strong solution to problem \eqref{iterationsystem} with the external terms $(g_*,h_*)$.
Furthermore, by the standard energy method, we obtain
\begin{align*}
& \frac12 \frac{\mathrm{d}}{\mathrm{d}t}\left( \|\nabla u\|^2+  \|\nabla v\|^2  + \alpha \sigma \|v-\overline{v} \|^2\right) +\|\nabla \mu\|^2 +
 \|\nabla \varphi\|^2 +\frac{\alpha}{2}\|\partial_t u\|^2 +\frac{\alpha}{2}\|\partial_t v\|^2 + \frac{\sigma}{2} \|\nabla v \|^2 \\
 &\quad \leq \frac{1}{2\alpha}\|g_*\|^2+ \frac{1}{2\alpha}\|h_*\|^2 +\alpha\sigma |\Omega|(\overline{v}-c)^2,
\end{align*}
and
 \begin{align*}
        & \frac{1}{2} \frac{\mathrm{d}}{\mathrm{d}t}\left( \| \nabla \mu \|^2 + \| \nabla \varphi\|^2 + \alpha \| \partial_t u\|^2 + \alpha \| \partial_t v\|^2+ 2\sigma \int_\Omega (v-c)\varphi \,\mathrm{d}x
        - \sigma\|\nabla v \|^2\right) \\
        &\qquad
          +   \|\nabla \partial_t u\|^2 +   \|\nabla \partial_t v\|^2    \\
        &\quad \le  C\left(\|\partial_t u\|^2+ \|\partial_t v\|^2+ \|g_*\|^2+\|h_*\|^2\right),
    \end{align*}
where
\begin{align*}
2\sigma \left|\int_\Omega (v-c)\varphi \,\mathrm{d}x\right|\leq \sigma\|\nabla v\|^2 + \alpha \sigma \|v-c\| \|\partial_t v\|  + \sigma\|v-c\|\|h_*\|.
\end{align*}
The above estimates enable us to conclude that $\mathcal{M}$ maps bounded set in $L^2(0,T; W)^2 \cap H^1(0,T;H)^2$ to bounded set in $L^2(0,T;H^3)^2\cap H^1(0,T;V)^2$.

Thanks to the nice property of $F_\delta$, with a refined argument, we can verify the existence of a strong solution to problem \eqref{eq:as1} on $[0,T]$ (see Definition \ref{ssad}) by using e.g., the Schauder principle. On the other hand, uniqueness of strong solutions follows from an easy application of the energy method. Since the argument is standard, we omit the details here.
\end{proof}

\subsection{Results for a nonlinear elliptic problem}

Let us consider the Neumann problem for an elliptic equation with singular nonlinearity:
\begin{equation}\label{eq:es}
\left\{
    \begin{aligned}
        &-\Delta u + \widehat{S}'(u) = f &&\quad  \mathrm{in} \ \Omega,\\
        & \partial_\mathbf{n} u = 0 &&\quad  \mathrm{on} \ \partial \Omega.
    \end{aligned}
\right.
\end{equation}
The following useful result can be found in \cite[Appendix]{GMT2019} (see also \cite{A2009,GGM2017}).

\bl \rm \label{nonlinearelliptic}
Let $\Omega$ be a bounded domain in $\mathbb{R}^d$  $(d\in\{2,3\})$  with smooth boundary $\partial \Omega$. Suppose that $\widehat{S}$ satisfies the assumption $\mathbf{(H1)}$.
For any $f \in H$, problem \eqref{eq:es} admits a unique solution $u \in W$ with   $\widehat{S}'(u) \in H$ such that it satisfies $-\Delta u + \widehat{S}'(u) = f $  almost everywhere in $\Omega$ and the boundary condition $\partial_\mathbf{n} u = 0$ almost everywhere on $\partial \Omega$.
Besides, we have the following estimate
\begin{equation}\label{eq:ellipticestimate}
    \| u\|_{H^2(\Omega)} + \|\widehat{S}'(u)\| \le C(1+\|f\|).
\end{equation}
If $f \in L^p(\Omega)$, $p\in [2,+\infty)$, it holds
    \begin{equation*}
        \|\widehat{S}'(u)\|_{L^p(\Omega)} \le \| f\|_{L^p(\Omega)}.
    \end{equation*}
In addition, if $f \in V$, we have
    \begin{equation*}
       \|\Delta u \|  \le \|\nabla u \|^{\frac{1}{2}}\|\nabla f\|^{\frac{1}{2}},
    \end{equation*}
and
\begin{equation*}
       \| u \|_{W^{2,p}(\Omega)} + \|\widehat{S}'(u)\|_{L^p(\Omega)} \le C(1+\| f\|_V),
    \end{equation*}
    for $p=6$ if $d=3$, and $p\in [2,+\infty)$ if $d=2$.
\el

\section*{Declarations}
\noindent \textbf{Acknowledgment.}  
The author would like to thank Professor Hao Wu for his encouragement and valuable suggestions. This work is partially supported by National Natural Science Foundation of China under Grant number 12071084 (PI: Hao Wu). \smallskip 

\noindent \textbf{Conflict of Interest.} 
The author declare that they have no conflict of interests. \smallskip

\noindent \textbf{Data Availability.} 
Data sharing not applicable to this article as no datasets were generated or analysed during the current study.


\end{document}